\documentclass[11pt]{article}
\newcommand		{\comment}[1]		{}

\usepackage{amsmath,amsthm} 			
\usepackage{mathtools}

\usepackage[left=2.54cm,right=2.54cm,top=2.54cm,bottom=2.54cm]{geometry}
\usepackage[nouppercase]{scrlayer-scrpage}
\usepackage{xspace}
\usepackage{enumitem}

\usepackage{avant} 				
\usepackage[sc,osf]{mathpazo} 			
	\AtBeginDocument{
		\DeclareSymbolFont{AMSb}{U}{msb}{m}{n}
		\DeclareSymbolFontAlphabet{\mathbb}{AMSb}
	}
\usepackage{extpfeil} 
\usepackage{mathabx}  
\usepackage{mathrsfs}	
\usepackage{fancyvrb}	

\usepackage{rotating,pdflscape}
\usepackage{afterpage}
\usepackage{graphicx}
\usepackage{float}
\usepackage{subcaption}
\usepackage{caption}
\newcommand{\mockalph}[1]{\!}

\usepackage[utf8]{inputenc}
\usepackage[T1]{fontenc}

\usepackage{tocloft}
\makeatletter
\renewcommand{\l@figure}{\@dottedtocline{1}{1em}{3.5em}}
\renewcommand{\l@table}{\@dottedtocline{2}{1em}{3.5em}}
\makeatother
\newcommand*{\noaddvspace}{\renewcommand*{\addvspace}[1]{}}
\addtocontents{lof}{\protect\noaddvspace}

\usepackage[hyphens]{url}		%
\usepackage{hyperref}	
\usepackage{cleveref} 	
\newcommand		{\myred}		{BrickRed}
\hypersetup{
			colorlinks		= true, 	
			urlcolor		= \myred, 
			linkcolor		= \myred, 	
			citecolor		= PineGreen 	
}
\newcommand		{\hyref}[1]		{\hyperref[#1]{\ref*{#1}}}

\newif\ifdebug
\ifdebug\usepackage{lineno}\linenumbers\else\fi

\usepackage{etoolbox}
\PassOptionsToPackage{dvipsnames,svgnames,x11names,table}{xcolor}
\newcommand		{\defd}[1]	{\textcolor{RoyalBlue}{\textbf{\textit{#1}}}}
\newcommand		{\defm}[1]	{\textcolor{RoyalBlue}{#1}}

\usepackage[all]{xypic}\xyoption{rotate}
	\newdir{ >}{{}*!/-7pt/\dir{>}}
	\newdir{i}{{}*!/-7pt/\dir{(}	}
\usepackage{tikz}
\usetikzlibrary{matrix,fit,backgrounds,calc,arrows} 
\tikzstyle{image}=[rectangle,fill=Red!20,inner sep=-2pt]
\tikzstyle{nonzero}=[rectangle,fill=Navy!20,inner sep=0pt]
\tikzstyle{nonzerosm}=[rectangle,fill=Navy!20,inner sep=-2pt]

\usepackage{thmtools}
\makeatletter
\let\c@figure\c@table
\let\c@equation\c@table
\makeatother

\numberwithin{table}{section}
\numberwithin{figure}{section}
\newtheorem{abcthm}[table]{Theorem}

\newtheorem{theorem}[table]{Theorem}

\newtheorem{proposition}[table]{Proposition}

\newtheorem{corollary}[table]{Corollary}

\newtheorem{lemma}[table]{Lemma}

\newtheorem{claim}[table]{Claim}

\theoremstyle{definition}

\newtheorem{definition}[table]{Definition}

\newtheorem{notation}[table]{Notation}

\newtheorem{observation}[table]{Observation}

\newtheorem{conjecture}[table]{Conjecture}

\theoremstyle{remark}
\newtheorem{fact}[table]{Fact}

\newtheorem{example}[table]{Example}

\newtheorem{exercise}[table]{Exercise}
\newtheorem{EG}[table]{Example}

\newtheorem{problem}[table]{Problem}
\newtheorem{histrmks}[table]{Historical remarks}
\newtheorem{remark}[table]{Remark}
\newtheorem{remarks}[table]{Remarks}

\theoremstyle{plain}
\newtheorem*{thm*}{Theorem}
\newtheorem*{theorem*}{Theorem}
\newtheorem*{prop*}{Proposition}
\newtheorem*{proposition*}{Proposition}
\newtheorem*{lemma*}{Lemma}
\newtheorem*{corollary*}{Corollary}
\newtheorem*{cor*}{Corollary}

\theoremstyle{definition}
\newtheorem*{definition*}{Definition}
\newtheorem*{defn*}{Definition}
\newtheorem*{QQ*}{Question}
\newtheorem*{obs*}{Observation}
\newtheorem*{notation*}{Notation}
\newtheorem*{discussion*}{Discussion}

\theoremstyle{remark}
\newtheorem*{rmk*}{Remark}
\newtheorem*{remark*}{Remark}
\newtheorem*{examples*}{Examples}
\newtheorem*{example*}{Example}
\newtheorem*{EG*}{Example}
\newtheorem*{EGs*}{Examples}
\newtheorem*{fact*}{Fact}
\newtheorem*{prob*}{Problem}

\newcommand{\bthm}{\begin{theorem}}
\newcommand{\ethm}{\end{theorem}}
\newcommand{\bprop}{\begin{proposition}}
\newcommand{\eprop}{\end{proposition}}
\newcommand{\bcor}{\begin{corollary}}
\newcommand{\ecor}{\end{corollary}}
\newcommand{\bconj}{\begin{conjecture}}
\newcommand{\econj}{\end{conjecture}}
\newcommand{\blem}{\begin{lemma}}
\newcommand{\elem}{\end{lemma}}
\newcommand{\bclm}{\begin{claim}}
\newcommand{\eclm}{\end{claim}}
\newcommand{\bpf}{\begin{proof}}
\newcommand{\epf}{\end{proof}}
\newcommand{\bdetails}{\begin{details}}
\newcommand{\edetails}{\end{details}}
\newcommand{\bdefi}{\begin{definition}}
\newcommand{\edefi}{\end{definition}}
\newcommand{\bdefn}{\begin{definition}}
\newcommand{\edefn}{\end{definition}}
\newcommand{\bex}{\begin{example}}
\newcommand{\eex}{\end{example}}
\newcommand{\bprob}{\begin{problem}}
\newcommand{\eprob}{\end{problem}}
\newcommand{\bob}{\begin{observation}}
\newcommand{\eob}{\end{observation}}
\newcommand{\bexer}{\begin{exercise}}
\newcommand{\eexer}{\end{exercise}}
\newcommand{\bexers}{\begin{exercises}}
\newcommand{\eexers}{\end{exercises}}
\newcommand{\brmk}{\begin{remark}}
\newcommand{\ermk}{\end{remark}}
\newcommand{\bhist}{\begin{histrmks}}
\newcommand{\ehist}{\end{histrmks}}
\newcommand{\brmks}{\begin{remarks}}
\newcommand{\ermks}{\end{remarks}}
\newcommand{\bverb}{\begin{verbatim}}
\newcommand{\everb}{\end{verbatim}}
\newcommand{\bntn}{\begin{notation}}
\newcommand{\entn}{\end{notation}}
\newcommand{\bfct}{\begin{fact}}
\newcommand{\efct}{\end{fact}}
\newcommand{\bfcts}{\begin{facts}}
\newcommand{\efcts}{\end{facts}}
\newcommand{\benum}{\begin{enumerate}}
\newcommand{\eenum}{\end{enumerate}}
\newcommand{\bitem}{\begin{itemize}}
\newcommand{\eitem}{\end{itemize}}

\tracingpatches
\makeatletter
\patchcmd{\@setref}{\bfseries ??}{\bfseries\color{red} FIX ME!}{}{}
\patchcmd{\@setcite}{\bfseries ?}{\bfseries\color{red} FIX ME!}{}{}
\patchcmd{\@setcref}         {??}{\color{red} FIX ME!}{}{}
\patchcmd{\@setcref}         {??}{\color{red} FIX ME!}{}{}
\patchcmd{\@setcrefrange}    {??}{\color{red} FIX ME!}{}{}
\patchcmd{\@setcrefrange}    {??}{\color{red} FIX ME!}{}{}
\patchcmd{\@setcrefrange}    {??}{\color{red} FIX ME!}{}{}
\patchcmd{\@setcrefrange}    {??}{\color{red} FIX ME!}{}{}
\patchcmd{\@setcrefrange}    {??}{\color{red} FIX ME!}{}{}
\patchcmd{\@setcrefrange}    {??}{\color{red} FIX ME!}{}{}
\patchcmd{\@setnamecref}     {??}{\color{red} FIX ME!}{}{}
\patchcmd{\@setnamecref}     {??}{\color{red} FIX ME!}{}{}
\patchcmd{\@setcpageref}     {??}{\color{red} FIX ME!}{}{}
\patchcmd{\@setcpageref}     {??}{\color{red} FIX ME!}{}{}
\patchcmd{\@setcpagerefrange}{??}{\color{red} FIX ME!}{}{}
\patchcmd{\@setcpagerefrange}{??}{\color{red} FIX ME!}{}{}
\patchcmd{\@setcpagerefrange}{??}{\color{red} FIX ME!}{}{}
\patchcmd{\@setcpagerefrange}{??}{\color{red} FIX ME!}{}{}
\patchcmd{\@setcpagerefrange}{??}{\color{red} FIX ME!}{}{}
\patchcmd{\@cref}            {??}{\color{red} FIX ME!}{}{}
\def\blx@citation@entry#1#2{%
  \blx@bibreq{#1}%
  \ifinlist{#1}{\blx@cites}
    {}
    {\listgadd{\blx@cites}{#1}%
     \blx@auxwrite\@mainaux{}{\string\abx@aux@cite{#1}}}%
  \ifinlistcs{#1}{blx@segm@\the\c@refsection @\the\c@refsegment}
    {}
    {\listcsgadd{blx@segm@\the\c@refsection @\the\c@refsegment}{#1}}%
  \blx@ifdata{#1}%
    {}%
    {\ifcsdef{blx@miss@\the\c@refsection}%
       {\ifinlistcs{#1}{blx@miss@\the\c@refsection}%
          {{\bfseries\color{red} cite:} }%
          {\blx@logreq@active{#2{#1}}}}%
       {\blx@logreq@active{#2{#1}}}}}
\def\blx@citeadd#1{%
  \ifcsdef{blx@keyalias@\the\c@refsection @#1}
    {\edef\blx@realkey{\csuse{blx@keyalias@\the\c@refsection @#1}}}
    {\def\blx@realkey{#1}}%
  \expandafter\blx@citation\expandafter{\blx@realkey}\blx@msg@cundefon
  \expandafter\blx@ifdata\expandafter{\blx@realkey}
    {\advance\blx@tempcnta\@ne
     \listeadd\blx@tempa{\blx@realkey}}
    {\ifnum\blx@tempcntb>\z@\multicitedelim\fi
     \expandafter\abx@missing\expandafter{\blx@realkey}%
     \advance\blx@tempcntb\@ne}}
\makeatother

\newcommand{\presectionskip}{-1.5\baselineskip}
\newcommand{\postsectionskip}{0.3\baselineskip}
\makeatletter
\usepackage{titlesec}
\renewcommand{\section}{\@startsection
  {chapter}{0}{0mm}
  {\presectionskip}
  {\postsectionskip}
  {\sffamily\huge}}
\renewcommand{\section}{\@startsection
  {section}{1}{0mm}
  {\presectionskip}
  {\postsectionskip}
  {\sffamily\LARGE}}
\renewcommand{\subsection}{\@startsection
  {subsection}{2}{0mm}
  {\presectionskip}
  {\postsectionskip}
  {\sffamily\Large}}
\renewcommand{\subsubsection}{\@startsection
  {subsubsection}{3}{0mm}
  {\presectionskip}
  {\postsectionskip}
  {\sffamily\normalsize}}
\renewcommand{\@seccntformat}[1]{\csname the#1\endcsname.\quad}
\newcommand\HUGE{\@setfontsize\Huge{30}{47}} 
\makeatother
\usepackage{titling}
  \titleformat{\chapter}[display]
  {\sffamily\Large}
  {Chapter {\HUGE\normalfont\thechapter}}    
  {1em}
  {\huge}
  \titleformat{\section}
  {\sffamily\Large}
  {\normalfont{\@setfontsize\Large{18}{47}\thesection.}}    
  {1em}
  {}
  \titleformat{\subsection}
  {\sffamily\large}
  {\normalfont{\Large\thesubsection}.}    
  {1em}
  {}
  \titleformat{\subsubsection}
  {\sffamily\normalsize}
  {\normalfont{\large\thesubsubsection}.}    
  {1em}
  {}

\renewcommand		{\SS}				{\textsection}
\newcommand		{\centeredalign}[1]		{\begin{center}\begin{minipage}{\linewidth}%
								\begin{align*}#1\end{align*}%
								\end{minipage}\end{center}%
							}

\newcommand		{\abs}[1]			{\left\vert #1 \right\vert}

\newcommand		{\rest}[2]			{#1 \upharpoonright #2}

\newcommand		{\quation}[1]			{\begin{equation} #1 \end{equation}}

\newcommand		{\eqn}[1]			{\begin{align*} #1 \end{align*}}

\def			\SPSB#1#2			{\rlap{\textsuperscript{#1}}\textsubscript{#2}}
\makeatletter
\def			\smallunderbrace#1		{\mathop{\vtop{\m@th\ialign{##\crcr
							   $\hfil\displaystyle{#1}\hfil$\crcr
							   \noalign{\kern3\p@\nointerlineskip}%
							   \tiny\upbracefill\crcr\noalign{\kern3\p@}}}}\limits}
\makeatother

\makeatletter
\newcommand{\subalign}[1]{%
  \vcenter{%
    \Let@ \restore@math@cr \default@tag
    \baselineskip\fontdimen10 \scriptfont\tw@
    \advance\baselineskip\fontdimen12 \scriptfont\tw@
    \lineskip\thr@@\fontdimen8 \scriptfont\thr@@
    \lineskiplimit\lineskip
    \ialign{\hfil$\m@th\scriptstyle##$&$\m@th\scriptstyle{}##$\crcr
      #1\crcr
    }%
  }
}
\makeatother

\makeatletter
\newcommand		{\oset}[3][0ex]			{%
								\raisebox{.175ex}{$%
								  \mathrel{\mathop{#3}\limits^{
								    \vbox to#1{\kern-2\ex@
								    \hbox{$\scriptstyle#2$}\vss}}}
								    $}%
							    }
\makeatother

\usepackage{faktor,xfrac}

\newcommand		{\smallhalf}{\textstyle\frac 1 2}

\newcommand		{\quotientbig}[2]	{\left.{\raisebox{.35em}{$#1$}} \!\!\Big/\!\! 
	{\raisebox{-.5em}{$#2$}}\right.}
\newcommand		{\quotientmed}[2]	{{\raisebox{.2em}{$#1$}}\  \!\!\big/\!\!\ 
										{\raisebox{-.2em}{$#2$}}}

\newcommand		{\dsp}		{\displaystyle}
\newcommand		{\nd}		{\noindent}

\newcommand		{\thesp}	{0.5ex}

\newcommand		{\mn}		{\mspace{-2mu}}
\newcommand		{\mnn}		{\mspace{-1mu}}

\newcommand		{\vpz}		{\vphantom{X^{5^2}z_3}}

\newcommand		{\ol}			{\overline}
\newcommand		{\os}			{\overset}
\newcommand		{\us}			{\underset}
\newcommand		{\ub}[2]		{\underbrace{#2}_{#1}}
\newcommand		{\ul}			{\underline}
\newcommand		{\wh}			{\widehat}

\newcommand		{\wt}			{\widetilde}

\newcommand		{\seq}[2][n]	{#2_1,\ldots,#2_{#1}}

\newcommand		{\mr}			{\mathrm}
\newcommand		{\bb}			{\mathbb}

\newcommand		{\f}			{\mathfrak}

\renewcommand		{\a}		{\alpha}

\renewcommand		{\b}		{\beta}

\newcommand		{\g}		{\gamma}
\renewcommand		{\d}		{\delta}
\renewcommand		{\epsilon}	{\varepsilon}
\newcommand		{\e}		{\epsilon}
\newcommand		{\z}		{\zeta}
\newcommand		{\h}		{\eta}
\renewcommand		{\t}		{\theta}

\newcommand		{\vk}		{\varkappa}
\renewcommand		{\l}		{\lambda}

\newcommand		{\s}		{\sigma}

\newcommand		{\w}		{\omega}

\newcommand		{\G}		{\Gamma}

\DeclareSymbolFont{cmletters}{OT1}{cmr}{m}{n}
\DeclareMathSymbol{\Ups}{\mathalpha}{cmletters}{"7}
\renewcommand		{\Upsilon}	{\Ups}

\newcommand		{\ceq}		{\coloneqq}
\newcommand		{\eqc}		{\eqqcolon}

\newcommand		{\concat}	{^{\frown}}
\let\union\cup%
\renewcommand		{\cup}		{\mspace{-1mu}\smile\mspace{-1mu}}
\let\inter\cap%
\renewcommand		{\cap}		{\mspace{-1mu}\frown\mspace{-1mu}}

\makeatletter
\DeclarePairedDelimiterX
			{\pmodx}[1]	{(}{)}{{\operator@font mod}\mkern6mu#1}
					\renewcommand{\pmod}{%
					  \allowbreak
					  \if@display\mkern18mu\else\mkern8mu\fi
						  \pmodx
					}
\makeatother

\newcommand		{\ab}		{^{\mathrm{ab}}}
\renewcommand		{\th}		{^{\mathrm{th}}}

\renewcommand		{\:}		{\colon}
\renewcommand		{\-}		{^{-1}}

\renewcommand		{\o}		{\circ}

\ExplSyntaxOn
\NewDocumentEnvironment{adjunctions}{O{}}
{
	\cs_set_eq:cN {@arraycr} \farin_arraycr:
	\keys_set:nn { farin/adjunction } { #1 }
	\begin{array}
		{
			@{ \hspace { \dim_eval:n { \l_farin_left_shift_dim + \l_farin_padding_dim } } }
			r
			@{ {\farin_strut:} \l_farin_symbol_tl {} }
			l
			@{ \hspace { \dim_eval:n { \l_farin_right_shift_dim + \l_farin_padding_dim } } }
		}
	}
	{
	\end{array}
}
\keys_define:nn { farin/adjunction }
{
	leftshift       .dim_set:N = \l_farin_left_shift_dim,
	leftshift       .initial:n = 0pt,
	rightshift      .dim_set:N = \l_farin_right_shift_dim,
	rightshift      .initial:n = 0pt,
	padding         .dim_set:N = \l_farin_padding_dim,
	padding         .initial:n = 6pt,
	symbol          .tl_set:N  = \l_farin_symbol_tl,
	symbol          .initial:n = \longrightarrow,
	verticalspacing .dim_set:N  = \l_farin_vertspac_dim,
	verticalspacing .initial:n = {3pt},
}
\cs_new_protected:Npn \farin_strut:
{
	\vrule height \dim_eval:n { \ht\strutbox + 1.2\l_farin_vertspac_dim }
	depth  \dim_eval:n { \dp\strutbox + \l_farin_vertspac_dim }
	width 0pt
}
\makeatletter
\exp_args:NNo \cs_new:Npn \farin_arraycr:
{
	\@arraycr\hline
}
\makeatother
\ExplSyntaxOff

\newcommand		{\dual}		{\mn^\vee}
\renewcommand		{\.}		{\cdot}

\newcommand		{\x}		{\times}

\newcommand		{\semidirect}	{\rtimes}

\newcommand		{\oplushigher}	{\mathbin{\raisebox{.85pt}{$\displaystyle\oplus$}}}
\newcommand		{\otimeshigher}	{\mathbin{\raisebox{.85pt}{$\displaystyle\otimes$}}}
\DeclareMathOperator*	{\otimesvariable}{%
			\mathchoice {\raisebox{.85pt}{$\displaystyle\otimes$}}
						{\raisebox{.85pt}{$\otimes$}}
						{\raisebox{0.7pt}{$\scriptstyle\otimes$}}
						{\raisebox{0.2pt}{$\scriptscriptstyle\otimes$}}
						}
\newcommand		{\tensor}	{\otimesvariable}
\newcommand		{\xt}[3]	{{#2}\us{#1}\otimeshigher{#3}}
\newcommand		{\direct}	{\oplushigher}
\newcommand		{\ox}		{\tensor}

\newcommand		{\+}		{\direct}

\newcommand		{\Direct}	{\bigoplus}

\newcommand		{\limit}	{\varprojlim}
\newcommand		{\colim}	{\varinjlim}

\newcommand		{\exterior}	{\Lambda}
\newcommand		{\ext}		{\exterior}

\newcommand		{\col}		{^{0,\bul}}
\newcommand		{\row}		{^{\bl,0}}
\newcommand		{\bul}		{\bullet}

\newcommand		{\Ei}		{E_\infty}

\DeclareMathOperator	{\bl}		{bl }

\DeclareMathOperator	{\diag}		{diag}

\DeclareMathOperator	{\rk}		{rk }
\DeclareMathOperator	{\im}		{im }

\DeclareMathOperator	{\coker}	{coker }

\DeclareMathOperator	{\Tor}		{Tor}

\DeclareMathOperator	{\Stab}		{Stab }

\DeclareMathOperator	{\Ad}		{Ad }

\DeclareMathOperator	{\Aut}		{Aut }

\newcommand		{\SO}		{\mr{SO}}

\newcommand		{\U}		{\mr{U}}
\newcommand		{\SU}		{\mr{SU}}

\newcommand		{\Spin}		{\mr{Spin}}

\makeatletter
\newbox\xrat@below
\newbox\xrat@above
\newcommand		{\xrightarrowtail}[2][]	{%
						  \setbox\xrat@below=\hbox{\ensuremath{\scriptstyle #1}}%
						  \setbox\xrat@above=\hbox{\ensuremath{\scriptstyle #2}}%
				  \pgfmathsetlengthmacro{\xrat@len}{max(\wd\xrat@below,\wd\xrat@above)+.6em}%
  						\mathrel{\tikz [>->,baseline=-.55ex]
              					   \draw (0,0) -- node[below=-2pt] {\box\xrat@below}
                            					    node[above=-2pt] {\box\xrat@above}
                    						   (\xrat@len,0) ;}
						}
\newbox\xrat@below
\newbox\xrat@above
\renewcommand		{\xtwoheadrightarrow}[2][]{%
						  \setbox\xrat@below=\hbox{\ensuremath{\scriptstyle #1}}%
						  \setbox\xrat@above=\hbox{\ensuremath{\scriptstyle #2}}%
				  \pgfmathsetlengthmacro{\xrat@len}{max(\wd\xrat@below,\wd\xrat@above)+.6em}%
						 \mathrel{\tikz [->>,baseline=-.55ex]
					                 \draw (0,0) -- node[below=-2pt] {\box\xrat@below}
					                                node[above=-2pt] {\box\xrat@above}
						                       (\xrat@len,0) ;}
		       				}
\makeatother
\newcommand		{\xmono}	{\xrightarrowtail}
\newcommand		{\mono}		{\xmono{\phantom{\ \, }}}
\newcommand		{\longmono}	{\xmono[]{\ \ \ \ }}

\newcommand		{\xepi}		{\xtwoheadrightarrow}
\newcommand		{\epi}		{\xepi{\phantom{\ \, }}}
\newcommand		{\longepi}	{\xepi[]{\ \ \ \ }}

\newcommand		{\longto} 	{\longrightarrow}
\newcommand		{\lt}		{\longto}

\newcommand		{\from}		{\leftarrow}

\newcommand		{\lmt}		{\longmapsto}

\newcommand		{\simto}	{\xrightarrow{\sim}}

\newcommand		{\longsimto}	{\os\sim\longto}
\newcommand		{\isoto}	{\longsimto}

\newcommand		{\bij}		{\leftrightarrow}
\newcommand		{\longbij}	{\longleftrightarrow}

\newcommand		{\inc}		{\hookrightarrow}
\newcommand		{\xinc}		{\xhookrightarrow}
\newcommand		{\longinc}	{\xinc[]{\ \ \ \ }}
\DeclareRobustCommand	{\longincfrom}	{\text{\reflectbox{$\ \longinc$\ }}}

\newcommand		{\ideal}	{\unlhd}

\newcommand		{\hmt}		{\simeq}
\newcommand		{\iso}		{\cong}
\newcommand		{\homeo}	{\approx}

\newcommand		{\pt}		{\mr{pt}}

\newcommand		{\Z}		{\bb Z}
\newcommand		{\Q}		{\bb Q}
\newcommand		{\R}		{\bb R}
\newcommand		{\C}		{\bb C}

\newcommand 		{\HH}		{\mathbb{H}}

\DeclareMathOperator	{\id}		{id}

\newcommand		{\APL}		{A_{\mathrm{PL}}}

\newcommand		{\CGA}		{\textsc{cga}\xspace}

\newcommand		{\DGA}		{\textsc{dga}\xspace}

\newcommand		{\CDGA}		{\textsc{cdga}\xspace}

\newcommand{\Sp}{\mathsf{Sp}}

\newcommand		{\fk}		{{\mathfrak k}}

\newcommand		{\fs}		{{\mathfrak s}}
\newcommand		{\ft}		{{\mathfrak t}}

\newcommand		{\td}		{\tilde{d}}

\newcommand		{\tw}		{\widetilde{w}}

\newcommand		{\tA}		{\widetilde{A}}

\newcommand		{\tE}		{\widetilde{E}}

\newcommand		{\tG}		{\widetilde{G}}
\newcommand		{\tH}		{\widetilde{H}}
\newcommand		{\tK}		{\widetilde{K}}
\newcommand		{\tN}		{\widetilde{N}}

\newcommand		{\tS}		{\widetilde{S}}

\newcommand		{\wP}		{\widehat{P}}

\newcommand		{\cP}		{\widecheck{P}}

\renewcommand{\HH}{{ H \mn H}}

\renewcommand 		{\H}		{H^*}

\newcommand		{\hb}		{h^\bullet}

\newcommand		{\tEs}		{\widetilde{E}_6}

\newcommand		{\KGK}		{{}_K G_K}

\newcommand		{\SGS}		{{}_S G_S}
\newcommand		{\GK}		{$\smash{(G,K)}$\xspace}
\newcommand		{\GS}		{$\smash{(G,S)}$\xspace}

\newcommand		{\ewP}		{\ext\wP}

\newcommand		{\eP}		{\ext P}

\newcommand		{\ecP}		{\ext\cP}

\newcommand 		{\HG}		{{\H_G}}

\newcommand 		{\HT}		{{\H_T}}

\newcommand 		{\HK}		{{\H_K}}

\newcommand 		{\HS}		{{\H_S}}

\newcommand		{\hG}		{\hb(G)}
\newcommand		{\hGS}		{\hb(G/S)}
\newcommand		{\hGK}		{\hb(G/K)}

\newcommand		{\pN}		{N}
\newcommand		{\pNp} 		{N'}
\newcommand		{\apN}		{|\pN|}
\newcommand		{\apNp}		{|\pNp|}

\newcommand		{\SSS}		{Serre spectral sequence\xspace}

\newcommand		{\toral}	{toral\xspace}
\newcommand		{\circular}	{circular\xspace}

\newcommand		{\eqf}		{equivariantly formal\xspace}

\newcommand		{\isotf}	{isotropy-formal\xspace}

\numberwithin{equation}{section}

\let\lim\varprojlim

\newcommand{\Cech}{\widecheck{H}^*}
\newcommand{\wH}{\widecheck{H}}
\newcommand{\HKz}{\H_{K_0}}

\begin{document}
\title{\vspace{-1em}\huge Equivariant formality of isotropic torus actions}
\author{\Large Jeffrey D.~Carlson}
\maketitle

\begin{abstract}
\small{
It is an open question for which pairs $(G,K)$
of Lie groups $G$ and closed, connected subgroups $K$
the left action of $K$ on the homogeneous space $G/K$ 
is equivariantly formal. 
We arrive through a sequence of reductions
at the case
$G$ is compact and simply-connected
and $K$ is a torus.


To illustrate the feasibility of this approach, 
we classify all pairs $(G,S)$ such that $G$ is compact connected Lie 
and the embedded circular subgroup $S$ acts equivariantly formally on $G/S$.
In the process we provide what seems to be the first published proof of the 
(long known) structure of the cohomology rings $\H(G/S;\Q)$.
%
%
%
%
%
%
}
\end{abstract}

\section{Introduction}
A natural request of
a continuous group action $G \x X \lt X$
is that it be \emph{equivariantly formal},
meaning the fiber inclusion in the Borel fibration $X \to X_G \to BG$ 
induces a surjection $\H_G(X;\Q) \longepi \H(X;\Q)$
of Borel equivariant cohomology upon singular cohomology. 
While the term was only coined in 1997 
by Goresky, Kottwitz, and MacPherson~\cite{GKM1998},
the condition had already been alighted upon by Borel in 
Chapter XII of his Seminar~\cite{borel1960seminar}.
This condition makes available a comparatively tractable computation of $\HG(X;\Q)$
in terms of $G$-orbits of dimensions zero and one in the case there are only finitely many of each,
as well as, by definition, guaranteeing all classes of $\H(X;\Q)$ have equivariant extensions
in $\H_G(X;\Q)$, to which, for example, 
the localization theorems of Berline--Vergne/Atiyah--Bott~\cite{BV1982}\cite{AB1984} and 
Jeffrey--Kirwan~\cite{jeffreykirwan1995} can be applied.

As any orbit of a continuous action of a Lie group $G$ on a space $X$,
is a homogeneous space $G/\Stab_G(x)$, it is natural to
ask about equivariantly formal actions on such spaces.
The transitive $G$-action is only equivariantly formal 
if the isotropy group $K = \Stab_G(x)$ is of full rank, 
but some restriction of this action to a subgroup $H$ will
always be equivariantly formal.
For this to happen,  
$H$ cannot contain a strictly larger maximal torus than $K$ does,
so that the left action of $K$ is in some sense the ``largest'' action on 
$G/K$ which could conceivably be equivariantly formal. 
Assuming that $G$ is compact, it is known that the isotropy action of $K$ on $G/K$
is equivariantly formal if $K$ is of full rank in $G$~\cite[Proposition 1]{brion1998eqcohom},
if $\H(G;\Q) \lt \H(K;\Q)$ is surjective \cite[Thm.~A, Cor.~4.2]{shiga1996equivariant},
or if $(G,K)$ is a generalized symmetric pair with $K$ connected~%
\cite{goertschesnoshari2016},
but otherwise few examples of such actions seem to be known.
Nevertheless, 
the full-rank case has found wide application in symplectic geometry
	(see, e.g., the book of Ginzburg, Guillemin, and Karshon~\cite{GGK},
	in which equivariant cohomology is already mentioned 
	in the first page of the introduction 
	and occupies a thirty-one--page appendix).

We show this question can be reduced to the case $K$ is a torus.
For concision, if the isotropy action of $K$ on $G/K$ is equivariantly formal,
we call the pair $(G,K)$ \defd{isotropy-formal}.

\begin{restatable}{theorem}{SGS}\label{SGS}
If $G$ is a compact Lie group, $K$ a closed, connected subgroup, and $S$ any torus
maximal within $K$, 
then $(G,K)$ is isotropy-formal if and only if $(G,S)$ is.
\end{restatable}
This result reduces the question to a study of embeddings of tori in Lie groups, 
an already more feasible-looking endeavor.
Further, the question reduces to the case the commutator subgroup of $G$ 
is simply-connected.

\begin{restatable}{theorem}{simplyconnected}\label{simplyconnected}
Let $G$ be a compact, connected Lie group, $K$ a closed,
connected subgroup, 
$\tG$ a finite central covering of $G$, 
and $\tK_0$ the identity component of the preimage of $K$ in $\tG$.
Then $(G,K)$ is isotropy-formal if and only if $(\tG,\tK_0)$ is.
\end{restatable}
The question is largely dependent on the case $G$ where itself is simply-connected.
 
\begin{restatable}{theorem}{semisimplered}\label{semisimplered}
Let $G$ be a compact, connected Lie group, $G'$ its commutator subgroup,
$K$ a closed connected subgroup of $G$,
and $S$ a maximal torus in $K$.
Write $K' = K \cap G'$ and $S' = S \cap G'$ for the intersections with $G'$
and $K'_0$ and $S'_0$ for their respective identity components.
Then $(G,K)$ is isotropy-formal
if and only if 
\begin{enumerate}
\item the pair  $(G',K'_0)$ is isotropy-formal and
\item the inclusion $N_G(S) \longinc N_{G}(S'_0)$ induces an isomorphism of component groups.
\end{enumerate}
\end{restatable}
These reductions are proven in \Cref{sec:red},
with some additional partial reductions having to do with disconnected goups and general compact Hausdorff groups expounded in \Cref{sec:failures}.
The reductions achieved, in \Cref{algworks},
we are able to completely determine 
for $G$ any compact, connected Lie group and $S$ any circular subgroup
whether $(G,S)$ is isotropy-formal. 
A key condition turns out to be that
there exist an element of the $G$ conjugation by which acts as
$s \longmapsto s\-$ on $S$.
We say such an element \defd{reflects} $S$.


\begin{restatable}{proposition}{main}\label{theorem:main}
	Let $G$ be a compact, connected Lie group and $S$ a circular subgroup of $G$. 
	There are the following three mutually exclusive cases. 
	\begin{enumerate}
		\item The inclusion $S \longinc G$ surjects in cohomology and $S$ is not reflected in $G$.
		\item The inclusion $S \longinc G$ is trivial in cohomology and
		\begin{enumerate}
			\item[2a.]
			$S$ is reflected in $G$.
			\item[2b.]
			$S$ is not reflected in $G$. 
		\end{enumerate}
	\end{enumerate}
	Only in the last case is $(G,S)$ not isotropy-formal.
\end{restatable}

Reflected circles can classified entirely,
and from 
Propositions \ref{theorem:main},
\ref{theorem:Hsurjecttorus},
\ref{Atrivial}, and
\ref{reflprod},
one assembles the following result.

\begin{theorem}\label{algor}
Let $G$ be a compact, connected Lie group and $S$ a \circular subgroup of $G$.
If $S$ is not contained in the commutator subgroup $G'$ 
of $G$, then $(G,S)$ is isotropy-formal.
Otherwise,
we may assume by \Cref{simplyconnected}
that $G'$ is a product of simple Lie groups $K_j$. 
Pick for each a maximal torus containing the 
image $S_j$ of $S \inc G' \epi K_j$.
Then $(G,S)$ is isotropy-formal
if and only if for each $K_j$ there is 
an element of the Weyl group $W(K_j)$
reflecting $S_j$,
which is determined as laid out in \Cref{TABLE}.
%
\end{theorem}

\bigskip

\begin{table}[H]	
\caption{Reflected circles in simple Lie groups}\label{TABLE}
\begin{center}
\begin{tabular}{l | l}
Type of $K$		& 
The circle $S$ in $K$ is reflected \ldots\\
\hline
$A_n$ & 
when the exponent multiset $J$ satisfies $J = -J$. \\
$B_n$ & always. \\
$C_n$ & always. \\
$D_{2n}$ & always. \\
$D_{2n+1}$ & if $S$ is contained in a $D_{2n}$ subgroup.\\
$G_2$ & always.\\
$F_4$ & always.\\
$E_6$ & if $S$ is contained in a $D_4$ subgroup.\\
$E_7$ & always. \\
$E_8$ & always. \\
\end{tabular}
\end{center}
\end{table}

This table is compiled in \Cref{sec:refl}.

\begin{remarks}[Explanatory remarks on \Cref{TABLE}]
The notation $\defm J$ in the $A_n$ case is the multiset
of exponents $\seq a \in \Z$
such that the injection $S^1 \mono \U(1)^{\+ n} \inc \U(n)$
realizing a conjugate of $S$ as a circular subgroup of 
the block-diagonal maximal torus of $\U(n)$ 
is given by $z \longmapsto \diag(z^{a_1},\ldots,z^{a_n})$.
We write $\defm{-J}$ for the multiset $\{-a_j\}_{1 \leq j \leq n}$
whose entries are the opposites of those of $J$; 
that is to say, for each $a \in \Z$, the element $-a$ occurs in $-J$ 
with the same multiplicity that $a$ occurs in $J$.
For example, $[-1 \ 0 \ 1] \in \Z^3$ meets the condition $J = -J$ 
and $[ 2\ 1\ -3]$ does not. See \Cref{reflSU}.

In the $D_{2n+1}$ case, 
$S$ is contained in a $D_{2n}$ subgroup just if it
is conjugate into a subtorus $T^{2n} \x \{1\}$
of the standard maximal torus $T^{2n+1}$
whose Lie algebra is the block-diagonal
subspace $\f{so}(2)^{\+ 2n+1}$ of $\f{so}(4n+2)$.
See \Cref{reflSpin}.

The condition that a circle in $E_6$ be contained in a $D_4$ subgroup
manifests, within a given maximal torus $T^6$ of $E_6$, 
in a more intricate fashion. 
Precise statements are \Cref{thm:XLV} and \Cref{rmk:E6}.
\end{remarks}

As an example of \Cref{algor},
we can recover Shiga's 
characterization~\cite[Prop.~4.3]{shiga1996equivariant}
of circles in the unitary group yielding isotropy-formality.

\begin{EG}
	If $S$ is a circle in the unitary group $\U(n)$, 
	then $\big(\mspace{-1mu}\U(n),S\big)$ is or is not isotropy-formal as indicated in \Cref{U(n)}.
	\begin{table}[H]
		\caption{The classification for circles in $\U(n)$}\label{U(n)}
		\begin{center}
			\begin{tabular}{ l|r }
				\mbox{Embedding of }$S$		
				& Is $\big(\U(n),S\big)$ isotropy-formal? \\ \hline
				$S\not\leq \SU(n)$ 					& Yes \\     \hline   
				$S\leq \SU(n) \mbox{ and }\  J = -J$ 		& Yes \\   
				$S\leq \SU(n) \mbox{ and }\  J \neq -J$ 	& No \\  
			\end{tabular}
		\end{center}
	\end{table}
\end{EG}

\bcor[anonymous referee]
Let $G$ be a compact, connected Lie group and 
$K$ a subgroup isomorphic to $\SO(3)$ or $\SU(2)$.
Then $(G,K)$ is isotropy-formal.
\ecor
\bpf
This follows from \Cref{algor} 
because the maximal torus $S^1$ of $K$
is contained in the commutator 
subgroup $G'$ of $G$ 
and is already reflected in $K$ 
and hence \emph{a fortiori} in $G$.
\epf
\bpf[Alternate proof]
Koszul~\cite[2.2${}^{\mathrm o}$]{koszul1947homologie} and 
Stiefel (unpublished)
showed $\H G \lt \H K$ is always surjective in this case
	(Samelson~\cite{samelson3formCR} derives this from the fact 
	the Cartan $3$-form 
	given at the identity by $(u,v,w) \lmt B\big(u,[v,w]\big)$
	is natural up to a scalar factor)
so it follows \cite[Thm.~A, Cor.~4.2]{shiga1996equivariant} 
that $(G,K)$ is isotropy-formal.
\epf


A crucial step of in obtaining the key \Cref{theorem:main} is the following
structure theorem for $\H(G/S)$,
which turns out to mildly extend 
a result which can be pieced together from two \emph{Comptes Rendus} 
notes of Leray and Koszul,
a complete proof of which seems never to have been published.
In case the result may be of independent interest, 
we take the opportunity to provide a proof in \Cref{sec:nosurj}.

\begin{restatable}{theorem}{thmHGS}\label{thmHGS}
	Let $G$ be a compact, connected Lie group and $S$ a circular subgroup. 
	\benum
	\item If $H^1 G \longto H^1 S$ is surjective, 
	then $\H(G/S) \lt \H G$
	is injective and its image is 
	the exterior algebra $\ewP$ on the 
	intersection $\wP$ of $\ker\mspace{-2mu}\big(\mspace{-1mu}\H G \to \H S)$
	with the graded vector space $P$ of primtive elements of 
	the exterior Hopf algebra $\H G = \eP$.
	Noncanonically, there is a $z_1 \in H^1 G$
	whose image spans $H^1 S$ and
	\[
		\H(G/S) = \ewP \, \iso \, \quotientmed{\H G}{(z_1)}.
	\]
	
	\item  If $H^1 G \longto H^1S$ is zero, 
	then the image of $\H(G/S) \lt \H G$
	is the exterior algebra on a codimension-one 
	subspace $\wP$ of $P$
	and $P / \wP \iso \Q z_3$ is graded in degree $3$.
	The image of $\HS \lt \H(G/S)$ is the subalgebra
	$\Q[s]/(s^2)$ generated by a nonzero $s \in H^2(G/S)$,
	and there are noncanonical isomorphisms
	\[
	\H(G/S) \,\iso \, \ewP \, \tensor\, \frac{\Q[s]}{(s^2)}
			\,\iso\, \frac{\H G}{(z_3)} \, \tensor\, \frac{\Q[s]}{(s^2)}.
	\]
	\eenum
\end{restatable}

\nd\emph{Acknowledgments.} 
The author thanks 
	Robert Bryant,
	Omar Antol{\'i}n Camarena,
	Jason DeVito,
	Chi-Kwong Fok,
	Oliver Goertsches,
	Fulton Gonzalez, 
	Jim Humphries, 
	Yael Karshon,
	Marek Mitros,
	George McNinch,  
	Jay Taylor, and
	Mathew Wolak
	for useful conversations,
Christopher O'Donnell for programming help,
Christopher Allday, Andreas Arvanitoyeorgos, 
and the referee for editorial advice, 
Oliver Goertsches and Sam Haghshenas Noshari for three figures, 
the referee for painstakingly careful reading and helpful suggestions,
and the National Center for Theoretical Sciences in Taipei
for its hospitality during part of this research.
He would finally and especially like to thank his advisor, Loring W.~Tu,
for patiently listening to and critiquing these ideas over a period of
many months. It is largely the result of his critical interventions
that this paper ever achieved coherence.

\section{Background}

%

Associated to a continuous action of a topological group $K$ on a space $X$, 
~\cite[IV.3.3, p.~53]{borel1960seminar} 
is the \emph{(Borel)} \defd{equivariant cohomology $\HK(X)$},
the rational singular cohomology $\H(X_K;\Q)$
of the \defd{homotopy quotient}~\cite[Def.~IV.3.1, p.~52]{borel1960seminar} 
(or \emph{Borel construction})
\[
	\defm{{}_K X} = 
	\defm{X_K} \ceq 
	\frac{EK \x X}{(ek,x) \sim (e,kx)},
\] 
where $EK \to BK$ is a universal principal $K$-bundle.
Until the last appendix, 
all cohomology will be singular cohomology with \textbf{rational coefficients},
which will henceforth be suppressed in the notation.
We write $\defm{\HK}$ for the coefficient ring $\H(BK) = \HK(\pt)$. 
Associated to the homotopy quotient is a fiber bundle $X \to {}_S X \to BS$,
the \defd{Borel fibration}.
As noted in the introduction,
an action of a topological group $S$ on a space $X$ 
is said to be \defd{equivariantly formal}
if the fiber inclusion $X \longinc {}_S X$ 
in this fibration surjects in cohomology.\footnote{\ %
	Dating back to Hans Samelson's~\cite{samelson1941samelson} 
	\emph{nicht homolog 0} and ``$\not\sim 0$'',
	a space $F$	has been said to be \emph{(totally) nonhomologous to zero} 
	in a superspace $E$ if its inclusion induces an injection $H_* F \lt H_* E$.
	The inclusion has also been said to be 
	\emph{(totally) noncohomologous to zero} in the same event, 
	and the condition is abbreviated variously 
	\emph{TNHZ}, \emph{TNCZ}, and \emph{n.c.z.},
	notwithstanding the fact the map in cohomology is only injective 
	if it is an isomorphism.
	In the present work we maintain a respectful distance from this terminology.
	}
This condition is equivalent to the spectral sequence 
of this bundle collapsing at the $E_2$ page~\cite[Lem.~C.24, p.~208]{GGK}. 
Given a Lie group $G$ and closed subgroup $K$, 
we refer to the natural left $K$-action 
on the homogeneous space $G/K$ of left cosets
as the \defd{isotropy action}.
For brevity, when the isotropy action of $K$ on $G/K$ is equivariantly formal
we call the pair $(G,K)$ \defd{isotropy-formal}.

Given a Lie group $G$, 
we write 
$\defm{Z(G)}$ for its center,
$\defm{G'}$ for its commutator subgroup,
$\defm{G\ab} \ceq G/G'$ for its abelianization,
$\defm{W_G}$ for its Weyl group,
and $\defm{N_G(K)}$ and $\defm{Z_G(K)}$ 
respectively for the normalizer and the centralizer of a subgroup $K$ in $G$.
If $S$ is a torus in $G$,
we write $\defm N \ceq \pi_0 N_G(S)$ for the component group of its normalizer.
We write $\defm{\hb(X)} \ceq \sum_{n \geq 0} \dim_\Q H^n X$
for the total Betti number, and
denote subgroup containment by ``\defd{$\leq$}'', 
isomorphism ``\defd{$\iso$}'',
homotopy equivalence ``\defd{$\hmt$}'', 
and homeomorphism ``\defd{$\homeo$}''.

%

\comment{
It is easiest to lift the discussion of reflectibility in a semisimple group $K$ 
to its universal cover $\tK$.
If $S$ is a torus in $K$,
then since the fiber $F$ is discrete and central, 
the preimage of $S$ in $\tK$
is a compact abelian Lie group of equal dimension,
so the identity component $\tS$ of this preimage 
is another torus of that dimension. 
The torus $S$ is reflected in $K$ just if $\tS$ is reflected in $\tK$.

The case of simply-connected groups 
immediately reduces to that of simple groups.
Finally, we analyze the situation for simple groups. 
The work is mostly in the case $K = E_6$.

This last will be a corollary of the following proposition.
}

\subsection{Earlier work}
As noted in the introduction,
the question we are interested in could be asked in the late 1950s
but only received a name in the 1990s.
As of the beginning of this work, there were only the three known classes
of cases in the introduction and the following general results 
of Shiga and Takahashi.

\begin{theorem}[Hiroo Shiga \cite{shiga1996equivariant}]
Let $G$ be a compact Lie group, $K$ a closed, connected subgroup, and $N_G(K)$ the normalizer.
If $(G,K)$ is a Cartan pair and the map 
$\H(G/K)^{N_G(K)} \inc \H(G/K) \to \H(G)$
induced by $G \longepi G/K$ is injective, then 
$K$ acts equivariantly formally on $G/K$. 
\end{theorem}

The notion of \emph{Cartan pair}~\cite[(3) on p. 70]{cartan1950transgression} here
is not the notion due to {\'E}lie Cartan describing symmetric spaces,
but an algebraic condition on the (Henri) Cartan model 
for $G/K$ described in \Cref{sec:nosurj} 
which amounts to the space $G/K$ being formal in the sense of rational homotopy theory.
Visually, it corresponds to the tensor factorization 
$E_2 = E_2\row \ox E_2\col$
in the Serre spectral sequence of the
Borel fibration $G \to {}_K G \to BK$
persisting to the $\Ei$ page.
Shiga's theorem can be equivalently restated as follows.

\begin{proposition}[Shiga]
Let $G$ be a compact Lie group, 
$K$ a closed, connected subgroup,
and $N_G(K)$ the normalizer.
If $(G,K)$ is a Cartan pair and the map $\H_G \longto (\HK)^{N_G(K)}$ 
is surjective, then $K$ acts equivariantly formally on $G/K$. 
\end{proposition}
%

%

The result also has a partial converse. 
In a later-written but earlier-published technical report~\cite{shigatakahashi1995},
Shiga and Hideo Takahashi prove a partial converse.

\begin{theorem}[Shiga--Takahashi]\label{thm:ST}
Let $G$ be a compact group, 
$S$ a \toral subgroup,
and $N_G(S)$ the normalizer.
Suppose that $S$ contains regular elements of $G$ and $(G,S)$ is a Cartan pair. 
Then $S$ acts equivariantly formally on $G/S$
if and only if and the map $\H_G \longto (\HS)^{N_G(S)}$ is surjective. 
\end{theorem}

In work with Chi-Kwong Fok~\cite{carlsonfok2018}, 
we show that if \GK is isotropy-formal, then $G/K$ must be formal,
so the ``Cartan pair'' hypothesis is redudant.
The hypothesis on regular elements is also unnecessary,
and in later work~\cite[Cor.~5.17]{carlson2022k}, 
we show that $S$ can also be replaced by any closed, connected subgroup $K$
in the result. 
Although we do not need it in what follows, 
we state the strong version for here for reference.

\begin{theorem}
	Let $G$ be a Lie group, 
	$K$ a closed, connected subgroup,
	and $N_G(K)$ the normalizer. 
	Then \GK is \isotf if and only if $G/K$ is formal and
	$\H_G \longto (\HK)^{N_G(K)}$ is surjective.
\end{theorem}

Our trichotomy \Cref{theorem:main} about the case $K \iso S^1$ 
can actually be refactored through 
the Shiga--Takahashi result.
Noting that the regular element condition is unneeded,
and that $G/S$ is always formal for $S$ a circle 
by the classical results of \Cref{sec:nosurj}, the Shiga--Takahashi 
theorem \ref{thm:ST} 
reduces isotropy-formality of \GS
to study of the map $\HG \lt \HS$.
In this language,
\Cref{theorem:main}
can be reproven as follows:
one has $N = \pi_0 N_G(S)$ either trivial or $\{\pm 1\}$.
If it is trivial, then isotropy-formality is just that
$\HG \longto \HS$ is surjective,
which happens if and only if $\H(G) \longto \H(S)$ surjects~\cite[$1^{\o}$, p.~69]{cartan1950transgression}%
\cite[Cor., p.~139]{borelthesis}.
Otherwise $N \iso \{\pm 1\}$,
meaning exactly that $S$ is reflected in $G$ (\Cref{card2}), 
and $N$ acts as $s \longmapsto \pm s$ on $\Q[s] \iso \H(BS)$,
so that $\H(BS)^N = \Q[s^2]$; 
then one proves \Cref{imageHBK} to see 
$\H(BG) \longto \Q[s^2]$ is always surjective.

The way this is presented in \Cref{sec:trichotomy}, 
we use a well-known fixed point criterion for equivariant
formality (\Cref{fpdim}) 
and a computation of the vector space dimension of the 
cohomology of the fixed point set due to Goertsches (\Cref{GNdim}).
Whether reasoning through a dimension count
or through \Cref{thm:ST}, 
one way or another the crux of it is understanding
the cohomology of the maps $S \to G \to G/S \to BS \to BG$.

\section{Reductions}\label{sec:red}

In this section we undertake a series of reductions that 
ultimately localizes most of the difficulty in determining which pairs \GK
are \isotf in the case where $G$ is semisimple and $K$ a torus.
Two further reductions, from disconnected to connected groups 
and from connected compact groups to Lie groups,
only go through partially and are sequestered in 
\Cref{sec:failures}.

\subsection{Compact total group}\label{sec:compact}

Let $G$ be connected pro-Lie group and $H$ a closed, connected subgroup. 
%
By the Cartan--Iwasawa--Malcev theorem,
there exists a maximal compact subgroup $\defm{K_H}$ of $H$, unique up to conjugacy
\cite[Cor.~12.77
]{hofmannmorrispro}, 
which is necessarily connected, 
such that there is a homeomorphism
$H \homeo K_H \x \R^{\kappa}$ for some cardinal $\kappa$~\cite[Cor.~12.82
]{hofmannmorrispro}.
Likewise $G$ contains a maximal compact subgroup $\defm{K_G}$, 
which after conjugation can be chosen  to contain $K_H$.
In case $G$ is a Lie group, at least, this yields a reduction result.

\bprop\label{pro-Lie}
Suppose $G$ is a connected Lie group 
and $H$ a connected, closed subgroup,
with respective compact, connected subgroups $K_G$ and $K_H$,
the one containing the other.
Then $(G,H)$ is isotropy-formal
if and only if $(K_G,K_H)$ is.
\eprop
\bpf
To identify the maps $\H_{K_H}(K_G/K_H) \lt \H(K_G/K_H)$ 
and $\H_H(G/H) \lt \H(G/H)$,
it will be enough to see that in the commutative diagram
\[
\xymatrix@C=1.25em@R=3.25em{
\ \ \, K_G/K_H  \ar@{}[r]^(.35){}="a"^(.825){}="b"\ar"a";"b"^(.35)\a \ar[d]	&
\ \ \, G/K_H	 \ar@{}[r]^(.325){}="a"^(.825){}="b"\ar"a";"b"^(.38)\g	
\ar[d]& 
\ G/H  			\ar[d]	\\
 {}_{K_H} \mn K_G/K_H  	\ar[r]_\b		& 
 {}_{K_H} \mn G/K_H  	\ar[r]_\d			&  
 {}_{H} G/H,
}
\]
the horizontal maps are homotopy equivalences. 
A left--$K_G$-equivariant deformation retraction 
of $G$ to $K_G$ induces deformation retractions 
from $G/K_H$ to $K_G/K_H$
and from ${}_{K_H}G/K_H$ to ${}_{K_H}K_G/K_H$. 
%
%
%
	The fibers of the bundles $\d$ and $\e$
%
%
are $H/K_H$ and $(H/K_H) \x (H/K_H)$ respectively,
both homeomorphic to Euclidean space,
and $G/K_H$ and $G/H$ have the homotopy type of a CW complex
so the long exact sequences of homotopy groups and Whitehead's theorem
show $\d$ and $\e$ are homotopy equivalences.
\epf
\brmk
This proof of \Cref{pro-Lie} depends only on homotopy equivalence, 
so the statement remains the same if $\H$ is replaced in the definition of isotropy-formality 
by any contravariant homotopy functor.
\ermk

\subsection{Toral isotropy}\label{sec:toral}

To reduce to toral isotropy actions, 
we require some well-known isomorphisms
and the rarely remarked fact these isomorphisms are \emph{natural}.

Let $\xi_0\colon E_0 \to B_0$ be a fibration with homotopy fiber $F$
such that $\pi_1 B_0$ acts trivially on $\H F$.
We can form a slice category 
of fibrations over $\xi_0$ with homotopy fiber $F$
by taking as objects maps of fibrations $\xi \to \xi_0$ with homotopy fiber $F$
and as morphisms between $\xi' \to \xi_0$ and $\xi \to \xi_0$ 
maps of fibrations $\xi' \to \xi$
making the expected triangle commute up to homotopy.
Such a morphism entails a homotopy-commutative prism
\begin{equation}\label{overbundlemap}
\begin{gathered}
	\xymatrix@C=3em@R=3em{
		 E' 	\ar[r]_{h} 	\ar@/^.65pc/[rr]	\ar[d]^{\xi'}	
		&E 		\ar[r] 							\ar[d]^{\xi}
		&E_0 									\ar[d]^{\xi_0}
	   \\B' 
	   	\ar@{}[r]^(.105){}="a"^(1.955){}="b"
	   		\ar[r]^{\bar h} \ar@<-.25ex>@/_.65pc/"a";"b"			
		&B 	
		\ar@{}[r]^(.075){}="a"^(.9){}="b"	\ar"a";"b" 	
		&\,B_0.
	}
\end{gathered}
\end{equation}

\begin{lemma}[{\cite[Cor.~4.4, p.~88]{smith1967emss}%
}]\label{TNHZpullback}
Let $\xi_0\colon E_0 \to B_0$ be a fibration
such that the fiber inclusion $F \longinc E_0$ is surjective in cohomology
and $\pi_1 B$ acts trivially on $\H B$.
Then the fiber inclusion of any fibration $\xi\: E \to B$ over $\xi_0$ with
homotopy fiber $F$ is surjective in cohomology, and
there is an $\H E_0$-algebra isomorphism
\[
	\smash{\xt{\H B_0 }{\H B \!}{\!\H E_0 } \isoto \H E}
\] 
natural in the fibration $\xi$ over $\xi_0$.
\end{lemma}

We prove the result so as justify the naturality clause
we will need,
absent in the original.
\begin{proof}
Surjectivity of $\H E  \lt \H F$ is implied
by that of $\H E_0  \lt \H F$
since the fiber inclusion $F \lt E_0$ 
factors up to homotopy as $F \to E \to E_0$.
For the isomorphism, note that because of these surjections, 
the Serre spectral sequences of these fibrations collapse at the $E_2$ page.
Thus the ring map $\H B \ox_{\H B_0 } \H E_0  \lt \H E $
induced by the maps in the right square of (\ref{overbundlemap})
is equivalent on the level of $\H B_0 $-modules
to the canonical isomorphism
\[
	\xt{\H  B_0 }{\H B \!}{\!\big(\mn\H B_0 \ox \H F\big)} 
		\isoto 
	\H B \ox \H F,
\]
and so is itself an isomorphism.
%
For naturality, note that the ring map $h^*\colon \H E \longto \H E'$
is completely determined its restrictions to its tensor-factors $\H B$ and $\H E_0$
and that the commutative diagrams in cohomology
induced by the left square and top triangle of (\ref{overbundlemap}) 
respectively
imply these restrictions are  
$\bar h^* \colon \H B  \longto \H B'$ and $\id_{\H E_0 }$.
\end{proof}

The naturality in the following lemma 
follows from the standard proof 
by noting that a $K$-equivariant map $X \lt Y$
yields commutative squares 
\[
	\xymatrix@C=1.25em@R=2.75em{
		X/S \ar[r]\ar[d]	&X/N_K(S)\ar[r]\ar[d]&X/K\, \ar[d]\\
		Y/S \ar[r]			&Y/N_K(S)\ar[r]		&Y/K.
	}
\]

\begin{lemma}[{\cite[Lemma~III.1.1, p.~35]{hsiang}}]\label{Weyl}
Let $K$ be a compact, connected Lie group with maximal torus $S$
and Weyl group $W$,
and $X$ a free $K$-space.
Then there is a ring isomorphism, natural in $X$,
\[
	\H( X / K) \longsimto \H(X/S)^W.
\] 
\end{lemma}

%
%


\begin{lemma}[{\cite[Prop.~III.1, p.~38]{hsiang}}]\label{KSdecomp}
Let $K$ be a compact, connected Lie group with maximal torus $S$
and Weyl group $W$.
Then there are the following ring isomorphisms natural in $X$:
	\eqn{
		\HK(X)		 			&\longsimto \HS(X)^W,\\
		\xt{\HK}{\HS}{\HK(X)} 	&\longsimto \HS(X). 
	\qedhere}
\end{lemma}
\bpf
The first statement follows from
\Cref{Weyl} and the definitions.
The second follows from \Cref{TNHZpullback},
applied to the $K/S$-bundle $X_S \to X_K$
viewed as a bundle over $BS \to BK$;
alternately, 
as $W_K$ acts on $H^2_S$ as a reflection group, 
$\HS$ is a free module over $\HK \iso (\HS)^{W_K}$ 
by the Chevalley--Shephard--Todd theorem~\cite[p. 192]{kane}
and \Cref{freeECext} applies.
\end{proof}

\bcor\label{thm:invsurj}
Let $K$ be a compact, connected Lie group with maximal torus $S$
and $X \to Y$ a $K$-equivariant map.
Then $\vk_K\:\HK Y \lt \HK X$ is surjective if and only if 
$\vk_S\:\HS Y \lt \HS X$ is.
\ecor
\bpf
\Cref{KSdecomp} 
identifies $\vk_K$ with the 
map of Weyl-invariants $(\mspace{-1mu}\vk_S)^W$
and $\vk_S$ with the base extension $\id_{\HS} \ox_{\HK} \vk_S$.
If $\vk_S$ is surjective, 
then it follows by averaging that $\vk_K$ is as well,
since $\vk_S$ is $W$-equivariant
and $|W|$ is invertible in $\Q$.
On the other hand if $\vk_K$ is surjective,
then since the functor $\HS \ox_{\HK} -$ is right exact,
$\vk_K$ is surjective as well.
\epf

Finally, the following well-known lemma
follows from the preceding ones.

\begin{lemma}[{\cite[Prop.~C.26, p.~207]{GGK}}]\label{eqfKS}
If $K$ is a compact, connected Lie group and $S$ a maximal torus, and $K$ acts on a space $X$, 
then the action of $K$ is equivariantly formal if and only if the restricted action of $S$ is.
\end{lemma}
%

We can now prove the promised reduction.

\SGS*
\begin{proof}
By \Cref{eqfKS},
it is enough to show that $K$ acts equivariantly formally on $G/S$
if and only if it does on $G/K$.
To do so, we may apply \Cref{thm:invsurj}
to the map of right $K$-spaces $G \lt {}_K G$.
%
\end{proof}

\subsection{The dimension criterion}\label{sec:fp}

Equivariant formality can be reduced 
to a condition on total Betti numbers.

\begin{lemma}[{{\cite[Prop. XII.3.4, p. 164]{borel1960seminar}%
			\cite[Prop.~3.1, p.~81]{goertsches2012isotropy}}}]%
\label{fpdim}
An action of a torus $S$ on a topological space $X$
with finite total Betti number
is equivariantly formal if and only if $\hb(X) = \hb(X^S)$.
\end{lemma}

For later reference, note one inequality always holds:

\begin{lemma}[Borel,~{\cite[IV.5.5 (p. 62)]{borel1960seminar}}{\cite[Lem.~C.24]{GGK}}]%
\label{Borelrk}%
If a torus $S$ acts on a topological space $X$ with finite total Betti number,
then $\hb(X) \geq \hb(X^S)$.
\end{lemma}

Let $G$ be a compact Lie group and $S$ a torus in $G$.
As the fixed point set of the left action of $S$ on $G/S$
is the quotient group $N_G(S)/S$ 
of the normalizer,
we need to know when $\hGS = \hb\big(N_G(S)/S\big)$.
The latter number is easily 
expressed in terms of other quantities.
Recall that we denote by $Z_G(S)$ the centralizer of $S$ in $G$,
by ${W_K}$ the Weyl group of $K$,
and by $N$ the component group $\pi_0 N_G(S)$.

\begin{lemma}\label{pN}
Conjugation induces a natural injection 
$N \longmono \Aut S$. 
This induces homeomorphisms $N_G(S) \homeo \pN \x Z_G(S)$
and $(G/S)^S = N_G(S)/S \homeo \pN \x Z_G(S)/S$.
If $K$ is a closed, connected subgroup with maximal
torus $S$, there is a further homeomorphism
$(G/K)^S = N_G(S)K/K \homeo (\mspace{-1mu}N / W_K) \x Z_G(S)/S$.
Particularly, $\hb\big((G/K)^S\big) = 2^{\rk G - \rk K} \. |N|/|W_K|$.
\end{lemma}
\begin{proof}
The centralizer $Z_G(S)$ is connected since it is
the union of
the maximal supertori of $S$ in $G$.
As $Z_G(S)$ is the kernel of the continuous homomorphism 
$n \longmapsto (x \mapsto nxn\-)$
from $N_G(S)$ into the discrete group $\Aut S \iso \Aut {\Z^{\rk S}}$,
it is the identity component of $N_G(S)$.
Thus $N = N_G(S)/Z_G(S)$; 
the homeomorphisms follow because group components are homeomorphic.
As for $K$, 
one notes 
$\pi\: (G/S)^S \lt (G/K)^S$ 
can be equivalently written
as the surjection 
$N_G(S)/S \epi N_G(S)/N_K(S) = N_G(S)/\big(N_G(S) \cap K\big) \iso N_G(S)K/K$
with fibers $nN_K(S)/S = nW_K$.
It follows $(G/K)^S$ has $\apN/|W_K|$ components,
each homeomorphic to $Z_G(S)/\big(Z_G(S)\cap K\big) = Z_G(S)/S$.
But since $Z_G(S)/S$ is a compact, connected Lie group,
$\H\big(Z_G(S)/S\big)$ is an exterior algebra on $\rk G - \rk S$
generators by Hopf's theorem~\cite[Satz I, p. 23]{hopf1941hopf},
\epf

\begin{proposition}[Goertsches--Noshari~{\cite[Props. 2.1, 3.1]{goertschesnoshari2016}}]\label{GNdim}
Let $G$ be a compact, connected Lie group 
and $K$ a closed, connected subgroup.
Write $N = \pi_0 N_G(S)$.
Then $(G,K)$ is isotropy-formal
if and only if
\[
\hGK \leq  2^{\rk G -\rk K} \cdot \frac{\apN}{|W_K|}.
\]
\end{proposition}
\begin{proof}
	Let $S$ be a maximal torus of $K$.
	By \Cref{eqfKS},
	we may replace the $K$-action on $G/K$
	with the $S$-action.
	By Lemmas \ref{fpdim} and \ref{Borelrk}
	this action is \eqf if and only if
	$\hGK \leq \hb\big((G/K)^S\big)$,
	which is
	$2^{\rk G - \rk K} \. { \apN }/{|W_K|}$ by \Cref{pN}.
\end{proof}

\subsection{Torus--cross--simply-connected total group }\label{sec:simplyconnected}
The structure theorem for compact, connected Lie groups%
~\cite[Thm.~V.(8.1) \& Ex.~V.(8.7).6, p.~233, 238]{brockertomdieck}
states that each  
admits a finite central extension $\defm{p}\colon \defm{\tG} \longto G$ 
such that the abelianization exact sequence
$1 \to {\tG'} \to \wt G \to \wt G\ab \to 0$
splits on the level of topological groups.
If the kernel of $p$ is $\defm{F}$, 
we can write
$
G \iso 
\tG / F 
.
$ 
The total space $\tG$ (but not $p$ itself, if $A \neq 0$)
is uniquely determined up to isomorphism.

In determining which toral isotropy actions are equivariantly formal,
we will show we can replace $G$ with $\tG$ and 
the connected isotropy subgroup $K$ (which we can take to be a torus) 
with the identity component $\defm{\tK_0}$ 
of its preimage $\defm{\wt K} = p\-K = F\tK_0$.

\begin{proposition}\label{universalH}
These assumptions induce isomorphisms
${\H(G/K) \simto \H(\tG/\tK) \simto\H(\tG/\tK_0)}$.
\end{proposition}
This is a result of 
the following lemma 
and the homeomorphism $\tG/\tK \os\homeo\lt G/K$.
\begin{lemma}\label{pathtriv}
Let $\Gamma$ be a path-connected topological group, 
$F$ a central subgroup,
and $H$ another subgroup such that $FH/H$ is finite.
Then the covering $FH/H \to \Gamma/H \to \Gamma/FH$ 
induces an isomorphism
$\H(\Gamma/FH) \isoto \H(\Gamma/H)$.
\end{lemma}

\begin{proof}
As $F$ is central, 
the covering action of $f\mnn H \in f\mnn H/H$
is given by $\g H \cdot f\mnn H = \g f\mnn H = f\mnn\g H$,
left multiplication by $f$. 
But $\G$ being path-connected,
left translation by its any element is homotopic to the identity.
Thus~\cite[Prop.~3G.1]{hatcher}
\[
\H(\Gamma/FH) \iso \H(\Gamma/H)^{FH/H} = \H(\Gamma/H).\qedhere
\]
\end{proof}

The components of the normalizer are also preserved under this substitution.

\begin{proposition}\label{normcover}
Under the foregoing assumptions, the projection $p\colon \tG \longto G$
induces an isomorphism 
$N_{\wt G}(\wt K_0)/Z_{\wt G}(\wt K_0) \isoto N_G(K)/Z_G(K)$.
Particularly, if $S$ is a torus, 
$\pi_0 N_{\wt G}(\wt S_0) \eqc \defm {\wt N} \iso N = \pi_0 N_G(S)$.

\end{proposition}
\begin{proof}
As $p$ is a homomorphism, 
it sends $\smash{N_{\wt G}(\wt K) \longto N_G(K)}$.
We show this restriction is surjective and the preimage of $Z_G(K)$ is ${Z_{\tG}(\tK)}$.
For surjectivity, given $\tw \in p\-N_G(\wt K_0)$,
note $\tw 1 \tw\- = 1$ 
and $p(\tw \tK_0 \tw\-) = K$,
so ${\tw\tK_0\tw\- = \tK_0}$. 
For the preimage, note that if $\wt z \in p\-Z_G(K)$,
then $\wt z\wt k\wt z\-\wt k\- \in \ker p$ 
for each $\wt k \in \tK_0$;
since $\ker p$ is discrete and $\wt z 1 \wt z\- 1\- = 1$,
such a $\wt z$ centralizes $\tK_0$.
\end{proof}

These facts in hand, we conclude the proof of \Cref{simplyconnected}.

\simplyconnected*
\begin{proof}
Let $S$ be a maximal torus of $K$ and $\tS_0$ its connected lift in $\tK$.
We know from \Cref{GNdim} that $(G,K)$ is isotropy-formal if and only if 
\[
\hb(G/K) =  2^{\rk G - \rk S} |N| / |W_K|,
\]
and the analogous statement holds of $(\tG,\tK)$.
But evidently $\rk \tG = \rk G$ and $\rk \tK = \rk K$ 
and $W_K \iso W_{\tK}$;
from \Cref{universalH}, we know $\hb(\tG/\tK) = \hb(G/K)$;
and from \Cref{normcover}, we know $\tN \iso N$.
\end{proof}

In what follows we can therefore replace $G$ with a cover 
$\tG = \tG' \x \tG\ab$. 
For later, when we specialize to circles, we note the following corollary of \Cref{normcover}.

\begin{corollary}\label{reflcover}
Under these hypotheses, the torus $S$ is reflected in $G$ just if $\tS$ is reflected in $\tG$.
\end{corollary}

\subsection{Semisimple total group}\label{sec:semisimplered}
In this section, $G$ is a connected, compact Lie group,
$G'$ again its commutator subgroup,
and $G\ab$ its abelianization.
To separate out information about $G'$,
we will need another covering lemma similar in spirit to 
\Cref{pathtriv}.

%

\begin{lemma}\label{pathtriv2}
	Let $\G$ be a compact, connected Lie group, 
	$\Xi$ an abelian subgroup, and
	$S$ a torus in $\Xi$ such that $\Xi/S$ is finite.
	Then the covering $\Xi/S \to \G/S \to \G/\Xi$ 
	induces an isomorphism
	$\H(\G/S) \isoto \H(\G/\Xi)$.
\end{lemma}
\begin{proof}
	As $\Xi$ is abelian, it is contained in
	the centralizer $Z_\G(S)$, which is path-connected, so that its right action on $\G/S$
	is cohomologically trivial.
	Thus
	$
	\H(\Gamma/\Xi) \iso \H(\Gamma/S)^{\Xi/S} = \H(\Gamma/S).
	$
\end{proof}

Given subgroup $H$ of $G$,
the canonical short exact sequence
$G' \to G \to G\ab$
descends to a fiber bundle
$G'/(G' \inter H) \to G/H \to \coker(H \inc G \epi G\ab)$.

\begin{proposition}\label{tensor}
If $H$ is connected, this bundle 
has the cohomology of a trivial bundle.
\end{proposition}
\begin{proof}
%
Consider a finite central cover of the form $\wt G = \wt G' \x \wt G\ab$.
Let $\wt H$ be the full preimage of $H$ in $\wt G$ 
and $\wt H_0$ its identity component.
We will show 
$\wt G'/ (\wt G' \inter \wt H_0) \to \tG/\tH_0 \to \coker(\tH_0 \to \tG\ab)$ is a trivial bundle.
Then the K{\"u}nneth theorem will yield the desired
ring decomposition, for 
$\coker(H \to G\ab)$ and $\coker(\wt H_0 \to \tG\ab)$
are tori of the same dimension,
and $\H(G/H) \iso \H(\wt G/ \wt H_0)$
by \Cref{universalH},
while $\wt G' / (\wt G' \inter \wt H_0) \epi 
\wt G' / (\wt G' \inter \wt H)
\os\homeo\to G'/(G' \cap H)$
is a normal covering with 
covering action induced by right 
translation by central elements of $\wt G$,
so by \Cref{universalH}
again, 
$\H\big(G'/(G' \inter H)\big) \iso \H\big(\tG'/(\tG' \inter \tH_0)\big)$.

The short exact sequence 
$\im(\wt H_0 \to \tG\ab) \to \tG\ab \to \coker(\wt H_0 \to \tG\ab)$ 
of tori splits
on the level of topological groups.
Replacing $\tG\ab$ with the product in the expression $\tG = \tG' \x \tG\ab$,
the projection of $\tH_0$ to the cokernel component is trivial,
so $\wt G/\wt H_0$ is the direct product of ${\coker(\wt H_0 \to \tG\ab)}$ and
${\big(\wt G' \x \im(\tH_0 \to \tG\ab)\big)/\tH_0}$.
But the inclusion of $\wt G'/ (\wt G' \inter \wt H_0)$
into the latter is a continuous bijection of compact 
Hausdorff spaces, hence a homeomorphism.
\end{proof}


Now we can carry through the claimed near-reduction to the semisimple case.

\semisimplered*

\begin{proof}
Note that $S'_0$ is a maximal torus in $K'_0$,
so by \Cref{SGS} it is enough to show
$(G,S)$ is isotropy-formal if and only if $(G',S'_0)$
is and the condition on normalizers holds.

From the decomposition $G = G' \.Z(G)$,
it follows that
$N_G(\G) = N_{G'}(\G)\.Z(G)$ and 
$Z_G(\G) = Z_{G'}(\G)\.Z(G)$ 
for any subgroup $\G$, 
so that particularly 
$\pi_0 N_G(S'_0) \iso \pi_0 N_{G'}(S'_0) \eqc \defm{N'}$.
As $G'$ is normal in $G$, 
there is also a containment $N_G(S) \leq N_G(S'_0)$, 
and so an induced monomorphism $N \longmono N'$. 
Thus from \Cref{pathtriv2},
Borel's \cref{Borelrk} 
for the action of $S'_0$ on $G'/S'_0$ and \Cref{pN}, 
we see 
\begin{equation}\label{GNKS}
\hb(G'/S') = \hb(G'/S'_0) \geq 
\apNp\, 2^{\rk G' - \rk S'} \geq \apN\, 2^{\rk G' - \rk S'}.
\end{equation}
Because rank is additive under direct products,
\[
	\rk G - \rk S = 
	\big(\mn\rk G' + \rk Z(G)\big) - 
	\big(\mn \rk S'_0 + \rk \im(S \to G\ab)\big) =
	 \rk G' - \rk S'_0 + \rk \coker(S \to G\ab),
 \]
so multiplying (\ref{GNKS}) by $2^{\rk \coker(S \,\to\, G\ab)}$ yields,
by \Cref{tensor},
\begin{equation}\label{GNGS}
\hb(G/S) \geq \apNp \,2^{\rk G - \rk S} \geq \apN \,2^{\rk G - \rk S}.
\end{equation}
\Cref{GNdim} states that $(G,S)$ is isotropy-formal if and only if 
the inequalities (\ref{GNGS}) are in fact equalities, 
which is equivalent to (\ref{GNKS}) being equalities.
But by \Cref{GNdim} again,
this can only happen if $(G',S')$ is isotropy-formal and $N' \bij N$. 
\end{proof}

\begin{remark}
It can really happen that the inequality $\apN \leq \apNp$ is strict.
For instance, let 
$G = A \x G'$ for
$A = S^1$ and $G' = \SU(2)^{2}$, 
pick a circle $S^1$ in $\SU(2)$,
and let $T$ be the maximal torus $(S^1)^3$ of $G$
and
$
\smash{S = \big\{(z,w,zw\-) : z,w \in S^1\big\}}
$
a rank-two subtorus,
so that $S' = S'_0 = \smash{\big\{(1,w,w\-) : w \in S^1 \big\}}$.
Then $N' = W_{\SU(2)} \iso \Z/2$
but $\pN = 1$.
\end{remark}

\section{Circular isotropy}\label{algworks}

Now we can tackle the case $S$ is a circle.
This section demonstrates
the statements of \Cref{algor} and \Cref{TABLE} 
regarding equivariant formality of circle actions.

\subsection{The trichotomy}\label{sec:trichotomy}

Let $S \iso S^1$ be a circle subgroup of a compact, connected Lie group $G$.

\begin{proposition}\label{card2}
Then the cardinality of $\pi_0N_G(S)$ is $2$ if $S$ is reflected in $G$
and $1$ otherwise.
\end{proposition}
\bpf
This follows from \Cref{pN}
since $s \lmt s^{-1}$ is the only nontrivial continuous automorphism of $S^1$.
\epf

As $\wt{H}^* S^1 = H^1 S^1$ 
is one-dimensional, 
$\H G \lt \H S$ 
is either surjective or trivial. 

\begin{restatable}{proposition}{hsurjecttorus}\label{theorem:Hsurjecttorus}
The inclusion $S \longinc G$ is trivial in cohomology if and only if 
$S$ is contained in the commutator subgroup $G'$,
if and only if the map induced in $H^1$ by
$S \to G \to G\ab$ is trivial.
\end{restatable}

\begin{proof}\label{proof2.3}
Since $G'$ is the kernel of $G \lt G\ab \eqc A$,
it contains $S$ just if the composition $S \to G \to A$ is trivial.
If so, then of course the map $H^1 A \lt H^1 S$ is trivial.
If $S \to G \to A$ is nontrivial, 
then its image is a circle, 
so the induced map $\pi_1 S \lt \pi_1 A$
is nonzero and hence injective,
and so $H^1 A \lt H^1 S$ is surjective.
But this map is nontrivial just if $H^1 G \lt H^1 S$ is
since $H^1 A \lt H^1 G$ is an isomorphism,
as can be seen for example by using \Cref{universalH}
to pass to a finite cover $\tA \x \wt{G'}$ 
with $0 = \pi_1 \wt{G'} = H^1 \wt{G'} = H^1 G'$.
%
%
\end{proof}

%
%
%

We can now prove \Cref{theorem:main}.

\main*

\begin{proof}
Recall from \Cref{GNdim}
that $(G,S)$ is isotropy-formal just when
$\hb(G/S) \leq 
\apN\, 2^{\rk G -\rk S}$.
\Cref{thmHGS} imposes the constraint that 
$\smash{\hGS \in \mnn\big\{\frac 1 2 \hG, \hG\mnn\big\}}$
and 
\Cref{card2} that $\apN \in \{1,2\}$.
By \Cref{Borelrk}, it is impossible that both $\hGS = \frac 1 2 \hG$ and $\apN = 2$ simultaneously, 
so there are only the following three cases. 

\begin{enumerate}
	\item We have $\hGS = \frac 1 2 \hG$ and $\apN = 1$. 
	The action is equivariantly formal.
	\item We have $\hGS = \hG$, and
	\begin{enumerate}
		\item[2a.]
		$\apN = 2$. 	The action is equivariantly formal.
		\item[2b.]
		$\apN = 1$. 	The action is not equivariantly formal.\qedhere
	\end{enumerate}
\end{enumerate}
\end{proof}


It remains to determine when $|\pN| = 2$,
or in other words when $S$ is reflected in $G$.

\subsection{Classification of reflected circles}\label{sec:refl}

In this section, we determine what circular subgroups $S$ 
of compact, connected Lie groups $G$ are reflected.
First, we may assume $S$ lies in some fixed maximal torus of $T$,
since all maximal tori are conjugate 
and for any $g \in G$ one has $g N_G(S)g\- = N_G(gSg\-)$.
Further, we may represent reflections by Weyl group elements,
in that $\pN \leq \Aut S$ is naturally 
a quotient of $N_W(S) \leq W$.

\begin{lemma}[{\cite[Exercise IX.2.4, p.~391]{bourbakiLie7}%
\cite[Lemma~9.7, p.~20]{dwyer1998Lie}}]%
\label{invNGT}
Let $G$ be a compact, connected Lie group, $T$ a maximal torus, and $S$ a subtorus.
Any automorphism of $S$
induced from conjugation by an element of $N_G(S)$
is also induced by an element of $N_G(T) \cap N_G(S)$. 
\end{lemma}
%
%
Precisely, the inclusion $N_G(T) \cap N_G(S) \longinc N_G(S)$ 
induces maps  
\[\dsp
	\frac{N_G(T)}
		 {T}
			\longincfrom
	\frac{N_G(T) \cap N_G(S)} 
		 {T}
			\longepi
	\frac{N_G(T) \cap N_G(S)} 
		 {N_G(T) \cap Z_G(S)}
			\isoto
	\frac{N_G(S)}
		 {Z_G(S)}
	.
\]

\begin{corollary}\label{invWeyl}
A toral subgroup $S$ is reflected 
in a compact, connected Lie group $G$ if and only if 
some element of the Weyl group $W$ of $G$
acts as $s \longmapsto s\-$ on $S$.
\end{corollary}

From \Cref{reflcover}, 
we may replace $G$ 
with the product $A \x G'$ of a torus $A$ 
and a simply-connected Lie group $G'$,
but $A$ is irrelevant:

\begin{restatable}{proposition}{Atrivial}\label{Atrivial}
A toral subgroup $S$ is reflected in a compact, connected Lie group
$G$ 
if and only if it lies in and is reflected in the commutator subgroup $G'$. 
\end{restatable}
\bpf
Since the conjugation action of $A$ is trivial,
circles reflected by $G$ are already reflected by $G'$.
From Propositions \ref{theorem:Hsurjecttorus} and \ref{theorem:main}, 
we know any reflected $S$ in $G$ is contained in $G'$. 
\epf

Reflectibility of a torus in a semisimple group $H$ in turn
depends only on simple factors.


\begin{restatable}{proposition}{reflprod}\label{reflprod}
A \toral subgroup $S$ is reflected in a product $\prod H_j$ of Lie groups if and only if
each of its images $S_j$ under the factor projections to $H_j$ 
is reflected in $H_j$.
\end{restatable}

\begin{proof}
Since the homomorphisms $\prod H_j \longepi H_i$ 
preserve conjugacy and inversion,
if $(h_j) \in \prod H_j$ reflects $S$,
then $h_j$ reflects $S_j$.
On the other hand, if some $h_j \in H_j$ reflects each $S_j$, 
then $(h_j)$ reflects $\prod S_j$, which contains $S$.
\end{proof}

We can in fact restrict attention to a single element of the Weyl group.

\begin{proposition}\label{longestword}
A circular subgroup $S$ is reflected in a simple Lie group $H$ 
if and only if 
it is reflected by the longest word $w_0$ in the Weyl group $W$ of $H$.
\end{proposition}
\bpf
If $C$ is the closed Weyl chamber containing a given nonzero element 
$v \in \fs < \ft$,
then $-v$ lies the ``opposite'' closed Weyl chamber $-C$.
The orbit $W\.v$ meets $-C$ in exactly one point~\cite[Thm.~5.16]{adamsLiebook},
which must be $w_0 \.v$ since $w_0 \.C = -C$, 
so $\fs$ is reflected if and only if $w_0 \. v = - v$.
\epf

There is a representation-theoretic restatement of
the same condition.

\bcor\label{Fok}
A circular subgroup $S$ is reflected in a simple Lie group $H$ 
 if and only if 
the irreducible representation of $H$ 
determined by $S$ is self-dual.
\ecor
\bpf
	Identify $\ft$ with its dual $\ft\dual$ through the 
	$W$-invariant inner product
	and let $\l$ be an additive generator of
	the intersection of $\fs$ 
	with the weight lattice of $H$.
	Then $S$ is reflected if and only if $w_0\.\l = -\l$.
	But the dual to the irreducible representation with highest weight $\l$
	is that with highest weight $-w_0\.\l$.
\epf

\brmk
The original proof of the classification in \Cref{TABLE}
was unnecessarily intricate and involved a computer algebra 
verification at one point, and has been greatly simplified
through the arguments in \Cref{longestword} and \Cref{Fok}, 
due to Jay Taylor~\cite{MO:E6} and Chi-Kwong Fok (personal communication).
\ermk

To construct \Cref{TABLE} we
march case by case through the Killing--Cartan classification.

\begin{restatable}{proposition}{manyrefl}\label{manyrefl}
A maximal torus $T$ of a simple compact Lie group $G$ whose type is one of
\[
B_n,\quad
C_n,\quad
D_{2n},\ \ \ 
G_2,\quad
F_4,\quad
E_7,\quad
E_8
\] 
is reflected in $G$.
\end{restatable}
\bpf
The longest word $w_0$ acts as $-\mn\id$ on the vector space $\ft$ 
carrying the defining representation of $W$
precisely for Coxeter groups $W$ of these types%
~{\cite[Lem.~27-2, p.~283]{kane}} 
so $T$ is reflected by \Cref{longestword}.
Alternately, but relatedly, 
central involutions of a Weyl group $W$ reflect the maximal torus $T$%
~\cite[Thm.~1.8]{dwyer2001center}
and the center of $W$ is isomorphic to $\Z/2$ 
precisely for Coxeter groups $W$ of these types%
~\cite[Rmk.~1.9]{dwyer2001center}.
\epf

In the remaining cases, the longest word $w_0 \in W$ does not act as $-\mn\id$ on $\ft$,
so more work is required.

\bprop\label{graphauto}
A circular subgroup $S$ is reflected in a simple Lie group $H$ 
whose Weyl group has trivial center 
(\emph{viz.} one of type $A_n$, $D_{2n+1}$, or $E_6$)
if and only if there is some $w \in W$
such that $w\. \fs$ lies in the fixed point subalgebra $\ft^\t$ of the Cartan subalgebra
under an automorphism $\t \in \Aut \ft $ 
induced by a nontrivial diagram automorphism of the Dynkin diagram of $H$.
\eprop
\bpf
From \Cref{longestword} we know $S$ is reflected if and only if 
$\fs$ is fixed pointwise by the nontrivial automorphism $-w_0 \in \Aut \ft$.
As $w_0 = \Ad(n_0)$ for some $n_0 \in N_H(T)$,
we can extend $-w_0$ to $-\Ad(n_0) \in \Aut \fk$.
Outer automorphisms of $\f k$ 
are induced~\cite[Prop.~D.40, p.~498]{fultonharris} 
by graph automorphisms of the Dynkin diagram $\G$ of $H$
in the sense that $(\Aut \fk)/(\Ad H) \iso \Aut \G$.
Since $W$ acts simply transitively on Weyl chambers,
and $-w_0$ stabilizes 
but does not fix the positive closed Weyl chamber $C$,
the automorphism $-\Ad(n_0)$ of $\fk$ is not inner and 
hence its outer isomorphism class 
corresponds to a nontrivial automorphism $\t$ of $\Gamma$.
This means the induced $\t \in \Aut \ft$ is the restriction of 
$-\mspace{-1mu}\Ad(n_0 k) \in \Aut \fk$ for some $k \in N_H(T)$,
so that  $\t$ fixes $\Ad(k\-) \fs$.
\epf

It thus remains to find the fixed point subalgebras of 
nontrivial diagram automorphisms for Lie algebras of type $A_n$, $D_{2n+1}$, and $E_6$.
In all of these proofs, we use the fact that
the $W$-equivariant isomorphism $\ft\dual\isoto\ft$ 
induced by the invariant inner product
is also equivariant with respect to $\t = -w_0$,
and so identifies the fixed point subspaces $(\ft\dual\mn)^\t$ and $\ft^\t$.

\begin{figure}
	\caption{The graph involution of $A_n$}
	\label{fig:Afold}
	\centering
	\begin{minipage}[t]{.4\textwidth}	
		\begin{xy}
			0;<1pt,0pt>:
			(0, 42)="v_11",				c+<44pt, 0pt>="v_12",			c+<44pt, 0pt>="v_13",			c+<32pt, 0pt>="v_14",
			(0, 10)="v_21",				c+<44pt, 0pt>="v_22",			c+<44pt, 0pt>="v_23",			c+<32pt, 0pt>="v_24",
			"v_11"*\cir<4pt>{},			"v_12"*\cir<4pt>{},				"v_13"*\cir<4pt>{},				"v_14"*\cir<4pt>{},		c-<16pt, 0pt>*{\ldots},
			"v_21"*\cir<4pt>{},			"v_22"*\cir<4pt>{},				"v_23"*\cir<4pt>{},				"v_24"*\cir<4pt>{},		c-<16pt, 0pt>*{\ldots},
			{"v_11" + <0pt, 9pt>}*{\text{\scriptsize $\alpha_1$}},
			{"v_12" + <0pt, 9pt>}*{\text{\scriptsize $\alpha_2$}},
			{"v_13" + <0pt, 9pt>}*{\text{\scriptsize $\alpha_3$}},
			{"v_14" + <0pt, 9pt>}*{\text{\scriptsize $\alpha_\ell$}},
			{"v_21" - <0pt, 9pt>}*{\text{\scriptsize $\alpha_{2\ell}$}},
			{"v_22" - <0pt, 9pt>}*{\text{\scriptsize $\alpha_{2\ell -1}$}},
			{"v_23" - <0pt, 9pt>}*{\text{\scriptsize $\alpha_{2\ell -2}$}},
			{"v_24" - <0pt, 9pt>}*{\text{\scriptsize $\alpha_{\ell + 1}$}},
			{"v_11" + <4pt, 0pt>}; {"v_12" - <4pt, 0pt>} **@{-},
			{"v_12" + <4pt, 0pt>}; {"v_13" - <4pt, 0pt>} **@{-},
			{"v_21" + <4pt, 0pt>}; {"v_22" - <4pt, 0pt>} **@{-},
			{"v_22" + <4pt, 0pt>}; {"v_23" - <4pt, 0pt>} **@{-},
			{"v_14" - <0pt, 4pt>}; {"v_24" + <0pt, 4pt>} **@{-},
			\ar@<-4pt>@/_2pt/@{.>}{"v_11" - <0pt, 4pt>};{"v_21" + <0pt, 4pt>},
			\ar@<-4pt>@/_2pt/@{.>}{"v_21" + <0pt, 4pt>};{"v_11" - <0pt, 4pt>},
			\ar@<-4pt>@/_2pt/@{.>}{"v_12" - <0pt, 4pt>};{"v_22" + <0pt, 4pt>},
			\ar@<-4pt>@/_2pt/@{.>}{"v_22" + <0pt, 4pt>};{"v_12" - <0pt, 4pt>},
			\ar@<-4pt>@/_2pt/@{.>}{"v_13" - <0pt, 4pt>};{"v_23" + <0pt, 4pt>},
			\ar@<-4pt>@/_2pt/@{.>}{"v_23" + <0pt, 4pt>};{"v_13" - <0pt, 4pt>},
			\ar@<-4pt>@/_2pt/@{.>}{"v_14" - <0pt, 4pt>};{"v_24" + <0pt, 4pt>},
			\ar@<-4pt>@/_2pt/@{.>}{"v_24" + <0pt, 4pt>};{"v_14" - <0pt, 4pt>}
			\end{xy}
		\end{minipage}
		\begin{minipage}[t]{.1\textwidth}	
		\quad
		\end{minipage}
		\begin{minipage}[t]{.4\textwidth}		
			\begin{xy}
			0;<1pt,0pt>:
			(0, 42)="v_11",				c+<44pt, 0pt>="v_12",			c+<32pt, 0pt>="v_13",			c+<44pt, -16pt>="v_m",
			(0, 10)="v_21",				c+<44pt, 0pt>="v_22",			c+<32pt, 0pt>="v_23",			
			"v_11"*\cir<4pt>{},			"v_12"*\cir<4pt>{},				"v_13"*\cir<4pt>{},				c-<16pt, 0pt>*{\ldots},		"v_m"*\cir<4pt>{},
			"v_21"*\cir<4pt>{},			"v_22"*\cir<4pt>{},				"v_23"*\cir<4pt>{},				c-<16pt, 0pt>*{\ldots},
			{"v_11" + <0pt, 9pt>}*{\text{\scriptsize $\alpha_1$}},
			{"v_12" + <0pt, 9pt>}*{\text{\scriptsize $\alpha_2$}},
			{"v_13" + <0pt, 9pt>}*{\text{\scriptsize $\alpha_{\ell-1}$}},
			{"v_m" + <0pt, 9pt>}*{\text{\scriptsize $\alpha_\ell$}},
			{"v_21" - <0pt, 9pt>}*{\text{\scriptsize $\alpha_{2\ell-1}$}},
			{"v_22" - <0pt, 9pt>}*{\text{\scriptsize $\alpha_{2\ell -2}$}},
			{"v_23" - <0pt, 9pt>}*{\text{\scriptsize $\alpha_{\ell + 1}$}},
			{"v_11" + <4pt, 0pt>}; {"v_12" - <4pt, 0pt>} **@{-},
			{"v_13" + <4pt, 0pt>}; {"v_m" - <4pt, 0pt>} **@{-},
			{"v_23" + <4pt, 0pt>}; {"v_m" - <4pt, 0pt>} **@{-},
			{"v_21" + <4pt, 0pt>}; {"v_22" - <4pt, 0pt>} **@{-},
			\ar@<-4pt>@/_2pt/@{.>}{"v_11" - <0pt, 4pt>};{"v_21" + <0pt, 4pt>},
			\ar@<-4pt>@/_2pt/@{.>}{"v_21" + <0pt, 4pt>};{"v_11" - <0pt, 4pt>},
			\ar@<-4pt>@/_2pt/@{.>}{"v_12" - <0pt, 4pt>};{"v_22" + <0pt, 4pt>},
			\ar@<-4pt>@/_2pt/@{.>}{"v_22" + <0pt, 4pt>};{"v_12" - <0pt, 4pt>},
			\ar@<-4pt>@/_2pt/@{.>}{"v_13" - <0pt, 4pt>};{"v_23" + <0pt, 4pt>},
			\ar@<-4pt>@/_2pt/@{.>}{"v_23" + <0pt, 4pt>};{"v_13" - <0pt, 4pt>}		
			\end{xy}
		\end{minipage}
\end{figure}
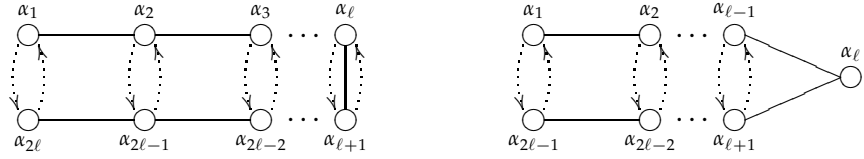

\begin{restatable}{proposition}{reflSUalg}\label{reflSUalg}
In a Lie algebra of type $A_n$, a point $v \in \ft\dual < \R^{n+1}$ 
of the dual Cartan algebra
is fixed by the automorphism $\t$ of \Cref{fig:Afold}
if and only if a permutation of the coordinates of $v$ yields $-v$.
\end{restatable}
\begin{proof}
The diagram automorphism $\t$ 
acts on simple roots of $A_n$ by exchanging $\a_j \longbij \a_{n-j}$.
The $\t$--fixed point subspace of $\ft\dual$ 
is spanned by the sums $\a_j + \a_{n-j}$
and so consists of those vectors 
$\sum c_j \a_j \in \ft\dual$
for which $c_j = c_{n-j}$.
The $\a_j$ are usually identified with $e_j - e_{j+1} \in \R^{n+1}$,
where $(e_\ell)_{1 \leq \ell \leq n+1}$ is the standard basis
and the resulting embedding $\ft\dual \longmono \R^{n+1}$ takes
\[
	\sum c_j \a_j 
		\lmt 
	\big[ c_1\ \ (c_2-c_1)\  \ \cdots\ \  (c_n-c_{n-1}) \ \ -c_n\big] 
		\eqc 
	\sum v_\ell e_\ell,
\]
translating the symmetry requirement $c_j = c_{n-j}$
to the antisymmetry condition $v_\ell = -v_{n+1-\ell}$.
\end{proof}

\vspace{-.25cm}

\begin{restatable}{corollary}{reflSU}\label{reflSU}
A circular subgroup $S$ is reflected in $\SU(n)$
if and only if 
the exponent multiset
$J$ of the inclusion of
any conjugate of $S$ into the standard maximal torus $T$ 
satisfies $J = -J$.
\end{restatable}

\bpf
Let $v$ span the tangent space $\fs < \ft$.
Recalling the Weyl group $W_{A_n} = S_{n+1}$ acts on $\R^{n+1}$ 
by permuting coordinates, 
by \Cref{reflSUalg}
a permutation of the entries of $v$ yields $-v$
just if some $w \in W_{A_n}$ sends $v$ into $\ft^\t$,
and by \Cref{graphauto}, $S$ is reflected just if this occurs.
\epf

\brmks
(\emph{a})
The root subsystems of $A_{2\ell}$ and $A_{2\ell-1}$ fixed by $\t$
are respectively of types $B_n$ and $C_n$,
corresponding to the inclusions $\SO(2\ell+1) \longinc \SU(2\ell+1)$
and $\Sp(n) \longinc \SU(2n)$
respectively induced by the ring injections $\R \longinc \C$ 
and $\HH \longmono \C^{2 \x 2}$.
These subgroups are fixed points of involutive automorphisms of $\SU(n)$
yielding the symmetric spaces
$\SU(n)/\SO(n)$ and $\SU(2n)/\Sp(n)$.
\medskip

\nd (\emph{b})
In terms of the self-duality criterion \Cref{Fok},
the representation $\tau$ of $S$ on $\C^n$ 
given by restricting the defining representation of $\SU(n)$  to $S$
is a direct sum $\smash{\Direct_{j=1}^n \rho^{\ox a_j}}$ 
of tensor powers of the defining representation 
$\rho \colon S^1 \longmono \Aut_\C \C$,
and the dual representation 
$\tau\dual = \smash{\Direct_{j=1}^n \rho^{\ox (-a_j)}}$,
will be isomorphic to $\tau$ just if $J = -J$. 

\ermks

%
%

\begin{figure}[H]
	\caption{The graph involution of $D_{2n + 1}$}
	\label{fig:Dfold}
	\centeredalign{
		\begin{xy}
			0;<1pt,0pt>:
			(0, 36)="v_1",		c+<44pt, 0pt>="v_2",	c+<32pt, 0pt>="v_3",	
			c+<44pt, +16pt>="v_u",		c-<0pt, 32pt>="v_d",
			"v_1"*\cir<4pt>{},		"v_2"*\cir<4pt>{},		
			"v_3"*\cir<4pt>{},		c-<16pt, 0pt>*{\ldots},		"v_u"*\cir<4pt>{},	"v_d"*\cir<4pt>{},
			{"v_1" + <0pt, 9pt>}*{\text{\scriptsize $\alpha_1$}},
			{"v_2" + <0pt, 9pt>}*{\text{\scriptsize $\alpha_2$}},
			{"v_3" + <0pt, 9pt>}*{\text{\scriptsize $\alpha_{2n-1}$}},
			{"v_u" + <0pt, 9pt>}*{\text{\scriptsize $\alpha_{2n}$}},
			{"v_d" - <0pt, 9pt>}*{\text{\scriptsize $\alpha_{2n + 1}$}},
			{"v_1" + <4pt, 0pt>};{"v_2" - <4pt, 0pt>} **@{-},
			{"v_3" + <4pt, 0pt>};{"v_u" - <4pt, 0pt>} **@{-},
			{"v_3" + <4pt, 0pt>};{"v_d" - <4pt, 0pt>} **@{-},
			\ar@<-4pt>@/_2pt/@{.>}{"v_u" - <0pt, 4pt>};{"v_d" + <0pt, 4pt>},
			\ar@<-4pt>@/_2pt/@{.>}{"v_d" + <0pt, 4pt>};{"v_u" - <0pt, 4pt>}
		\end{xy}
	}
\end{figure}
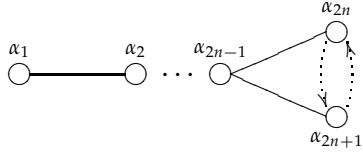

\begin{restatable}{proposition}{reflSpinalg}\label{reflSpinalg}
In a Lie algebra of type $D_{2n+1}$, 
a point $\l \in \ft\dual$ of the dual Cartan algebra
is fixed by an automorphism of the Dynkin diagram
if and only if the last coordinate of $\l$ is zero.
\end{restatable}
\begin{proof}
The nontrivial graph automorphism $\t$ of 
the Dynkin diagram of $D_{2n+1}$, shown in \Cref{fig:Dfold},
fixes all simple roots except $\a_{2n}$ and $\a_{2n+1}$, which it exchanges.
The fixed point subspace of $(\ft\dual\mn)^\t$ 
is spanned by $\{\a_j\}_{j < 2n} \union \{ \a_{2n} + \a_{2n+1} \} $.
The roots $\a_j$ for $j \leq 2n$ are usually identified 
with  $e_j - e_{j+1} \in \R^{2n+1}$
and $\a_{2n+1}$ with $e_n + e_{n+1}$,
where $(e_j)_{1 \leq j \leq n+1}$ is again the standard basis.
The image of the composite embedding $(\ft\dual)^\t \inc \ft\dual \to \R^{2n+1}$ 
is $\R^{2n} \x \{0\}$ 
since $\a_{2n} + \a_{2n+1} = 2e_{2n}$.
\end{proof}

\begin{restatable}{corollary}{reflSpin}\label{reflSpin}
A circular subgroup $S$ is reflected in $\Spin(4n+2)$
if and only if it is conjugate into
a $\Spin(4n)$ subgroup.
\end{restatable}

\bpf
Let $v$ span the tangent space $\fs < \ft = \R^{2n+1}$.
Recalling the Weyl group $W_{D_{2n+1}} = \{\pm 1\}^{2n} \semidirect S_{2n+1}$ 
acts on $\R^{2n+1}$ by permuting its coordinates and 
negating an even number of them, 
by \Cref{reflSpinalg}
some entry of $v$ is $0$
just if some $w \in W_{A_n}$ sends $v$ into $\ft^\t$,
and by \Cref{graphauto}, $S$ is reflected just if this occurs.
\epf

\brmk
The sublattice of a $D_{2n+1}$ lattice fixed by $\t$ is of type $B_{2n}$
and corresponds to a $\Spin(4n)$ subgroup of $\Spin(4n+2)$,
the fixed point set of an involutive automorphism of $\Spin(4n+2)$
yielding the symmetric space
$V_2(\R^{4n+2}) = \Spin(4n+2)/\Spin(4n) = \SO(4n+2)/\SO(4n)$.
%
%
\ermk

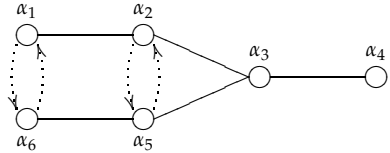
\begin{figure}[H]
		\caption{The graph involution of $E_6$}
		\label{fig:E6fold}
		\centeredalign{
			\begin{xy}
				0;<1pt,0pt>:
				(0, 42)="v_1",				c+<44pt, 0pt>="v_2",			c+<44pt, -16pt>="v_3",		c+<44pt, 0pt>="v_4",
				(0, 10)="v_6",				c+<44pt, 0pt>="v_5",			
				"v_1"*\cir<4pt>{},			"v_2"*\cir<4pt>{},				"v_3"*\cir<4pt>{},				"v_4"*\cir<4pt>{},		"v_5"*\cir<4pt>{},		"v_6"*\cir<4pt>{},   
				{"v_1" + <0pt, 9pt>}*{\text{\scriptsize $\alpha_1$}},
				{"v_2" + <0pt, 9pt>}*{\text{\scriptsize $\alpha_2$}},
				{"v_3" + <0pt, 9pt>}*{\text{\scriptsize $\alpha_3$}},
				{"v_4" + <0pt, 9pt>}*{\text{\scriptsize $\alpha_4$}},
				{"v_5" - <0pt, 9pt>}*{\text{\scriptsize $\alpha_5$}},
				{"v_6" - <0pt, 9pt>}*{\text{\scriptsize $\alpha_6$}},
				{"v_1" + <4pt, 0pt>}; {"v_2" - <4pt, 0pt>} **@{-},
				{"v_2" + <4pt, 0pt>}; {"v_3" - <4pt, 0pt>} **@{-},
				{"v_3" + <4pt, 0pt>}; {"v_4" - <4pt, 0pt>} **@{-},
				{"v_5" + <4pt, 0pt>}; {"v_3" - <4pt, 0pt>} **@{-},
				{"v_6" + <4pt, 0pt>}; {"v_5" - <4pt, 0pt>} **@{-},
				\ar@<-4pt>@/_2pt/@{.>}{"v_1" - <0pt, 4pt>};{"v_6" + <0pt, 4pt>},
				\ar@<-4pt>@/_2pt/@{.>}{"v_6" + <0pt, 4pt>};{"v_1" - <0pt, 4pt>},
				\ar@<-4pt>@/_2pt/@{.>}{"v_2" - <0pt, 4pt>};{"v_5" + <0pt, 4pt>},
				\ar@<-4pt>@/_2pt/@{.>}{"v_5" + <0pt, 4pt>};{"v_2" - <0pt, 4pt>}
			\end{xy}
		}
\end{figure}
\begin{restatable}{proposition}{reflE6alg}\label{reflE6alg}
In a Lie algebra of type $E_{6}$, a point $\l \in \ft\dual$ of the dual Cartan algebra
is fixed by the nontrivial automorphism of the Dynkin diagram
if and only if it lies in a certain $F_4$ sublattice.
\end{restatable}
\bpf
The fixed-point subspace $(\ft\dual\mn)^\t$
of the nontrivial automorphism $\t$ of the Dynkin diagram of $E_6$
depicted in \Cref{fig:E6fold}
is spanned by $\Delta = \{\a_1 + \a_6,\a_2+\a_5,\a_3,\a_4\}$.
By assumption, we have $\a_i \. \a_j = -2 |\a_i| |\a_j|$ for adjacent $\a_i,\a_j$ and $= 0$ otherwise,
so $\Delta$ is a simple root system of type $F_4$ with $a_1 + \a_6$ and $ \a_2+\a_5$ long 
and $\a_3$ and $\a_4$ short.
\epf

\begin{restatable}{proposition}{Esix}\label{Esix}
A \circular subgroup $S$ is reflected in $E_6$ or its universal cover $\tEs$ 
if and only if it is conjugate into a $\Spin(8)$ subgroup.
\end{restatable}
\bpf
It follows from \Cref{graphauto} and \Cref{reflE6alg}
that the tangent lines $\fs$ to reflected circles $S$ 
are precisely those sent into $\ft^\t$ by some $w \in W_{E_6}$.
As $(\ft\dual\mn)^\t$ is spanned by an $F_4$ sublattice of the $E_6$ root lattice,
its dual $\ft^\t$ is tangent to the maximal torus $T^4$ of an $F_4$ subgroup.
In the series of inclusions $\Spin(8) < F_4 < E_6$,
the first two share a maximal torus $T^4$, 
so $\ft^\t$ is actually tangent to the maximal torus of a $\Spin(8)$.
\epf

It may be of interest to count these four-dimensional tori.

\bprop\label{thm:XLV}
Within any given maximal torus $T^6$ of $E_6$ or $\tE_6$,
there are forty-five distinct Weyl-conjugate maximal tori 
$T^4$ of $\Spin(8)$ subgroups, all reflected.
\eprop
\bpf
The $\Spin(8)$ tangent to $T^4$ corresponds to a $D_4$
sublattice of $\ft$ spanning $\ft^\t$.
Within a set of positive roots for a root system of type $D_4$,
it is not hard to check
there are precisely three spanning sets of mutually orthogonal roots,
so the number of tori in question will be a third of the number 
of sets of four mutually orthogonal roots in the root system $\Phi(E_6)$.
Any given \emph{set} $\{\a,\b,\g,\d\}$ of four mutually orthogonal 
\emph{positive} roots in $\Phi(E_6)$ 
corresponds to $\abs{\{\pm 1\}^4 \semidirect S_4} = 384$ 
different mutually orthogonal \emph{ordered quadruples} of \emph{arbitrary} roots,
so the number of tori $T^4$ can be obtained by dividing the number of 
such quadruples by $384\.3 = 1152 = \big|W_{F_4}\big|$.
We will then be done if we can show $W_{E_6}$,
which is of cardinality $51,840 =  45\. 1152$,
acts simply transitively on mutually orthogonal
ordered quadruples $(\a,\b,\g,\d)$ in $\Phi(E_6)$.

For this, Carter observes~\cite[Lem.~11.(i), p.~14]{carter1972}
that $W_{E_6}$ acts transitively on roots $\a \in \Phi(E_6)$,
that $\Stab_{W_{E_6}} \a$ acts transitively on the $A_5$ subsystem
of roots $\b$ orthogonal to $\a$,
and that $\Stab_{W_{E_6}} (\a,\b)$ 
acts transitively on the $A_3$ subsystem 
the roots $\g$ orthogonal to both $\a$ and $\b$,
so that $W_{E_6}$ acts transitively
on mutually orthogonal ordered triples $(\a,\b,\g)$.
From there we may further see $\smash{\Stab_{W_{E_6}} (\a,\b,\g)}$
acts transitively on the $A_1$ subsystem $\{\pm \d\}$ 
of roots orthogonal to all of $\a,\b,\g$.
That the transitivity on quadruples is simple
follows, since $\big|\Phi(A_1)\big| = 2$, 
from repeated applications of the orbit--stabilizer theorem:
\[
	\ub{{|W_{E_6}|}}{51,840} 
		= 
	\ub{|\mn\Stab \a|}{720}\!\.\!\ub{|\Phi(E_6)|}{72} 
		= 
	\!\!\ub{|\mn\Stab (\a,\b)|}{24}	\!\!\!\.\,\!\! \ub{|\Phi(A_5)|}{30} \.	 	72 
		= 	
	\!\!\!\!\ub{|\mn\Stab (\a,\b,\g)|}{2}	\!\!\!\!\.\,\! \ub{|\Phi(A_3)|}{12} \.\; 30 \.	72.
	\qedhere
 \]
\epf

\brmk
If we view $T^4$ as the maximal torus of $F_4 < E_6$,
it follows from the equation $|W_{E_6}| = 45\.|W_{F_4}|$ 
that $W_{F_4}$ injects into $W_{E_6}$ as the normalizer of $T^4$.
The author is advised this result 
can be understood from Carter's book~\cite[Sec.~13.3]{carterbook}.
\ermk

\begin{remark}\label{rmk:E6}
A standard system of simple roots for $E_6$ in $\R^5 \x \R^3$ is given%
~
\cite[Planche~V, p.~260]{bourbakiLie4}
by
\[
\begin{array}{rr@{\hspace{\thesp}}r@{\hspace{\thesp}}r@{\hspace{\thesp}}r@{\hspace{\thesp}}
rrrl}	
		\defm\Delta	\ceq \big\{\mspace{.75mu} 
		\quad\ \ \defm\z &\ceq\  \smallhalf\mspace{1.5mu}  [1 & 1 & 1 & 1 & 1;& 1 & 1 & 1],\\
		-\defm{\g_{12}} &\ceq			\,					 - [1 & 1 & 0 & 0 & 0;& 0 & 0 & 0],\\
		\defm{\d_{12}} &\ceq			\,			\phantom{-}[1 &-1 & 0 & 0 & 0;& 0 & 0 & 0],\\
		\d_{23} &\ceq				\,			\phantom{-}[0 & 1 &-1 & 0 & 0;& 0 & 0 & 0],\\	
		\d_{34} &\ceq				\,			\phantom{-}[0 & 0 & 1 &-1 & 0;& 0 & 0 & 0],\\
		\d_{45} &\ceq			\,			\phantom{-}[0 & 0 & 0 & 1 &-1;& 0 & 0 & 0]\big\}.\\
\end{array}
\]
These roots span the six-dimensional subspace 
$
	\big(\mspace{-1mu}\R^5 \x \{\mspace{-1mu} 0\}^{\mn 3}\mspace{-1mu}\big)
		\; +\,\ 
	\R\cdot [1 \ \ 1 \ \ 1 \ \ 1 \ \ 1;\mspace{2mu}  1\ \ 1\ \ 1]
$
of $\R^8$ and one obtains a system $\Phi$ of 72 roots
obtained from permutation of the first five coordinates of 
\[
	\z,\quad
	\g_{12},\quad
	\d_{12},\quad
	\defm{\h_{12}} \ceq \z-\g_{12},\quad
	\defm{\e_1} \ceq \z - 2\g_{12} + 2\d_{12} + 3\d_{23} + 2 \d_{34}+ \d_{45}
			= \smallhalf\,[1\ {-1}\ {-1}\ -1\ -1;\ 1\ \ 1\ \ 1].		
\]
and multiplication by $\pm 1$.
We may choose the {positive} roots $\Phi^+$ to be the 36 in the union of 
the following 135 maximal mutually orthogonal sets: 
\eqn{
(60) &&
\{\e_a,\h_{ab},\g_{ac},\d_{de}\}, & 
	\qquad \mbox{ where }  \big|\{a,b,c,d,e\}\big| = 5 \mbox{ and } d<e,\\
(30) &&
\{\h_{ab},\h_{cd},\g_{ac},\g_{bd}\}, &
	\qquad \mbox{ where }  \big|\{a,b,c,d\}\big| = 4,\\
(15) &&
\{\h_{ab},\h_{cd},\d_{ab},\d_{cd}\}, &
	\qquad \mbox{ where }  \big|\{a,b,c,d\}\big| = 4 \mbox{ and } a<b \mbox{ and } c<d,\\
(15) &&
\{\g_{ab},\g_{cd},\d_{ab},\d_{cd}\}, &
	\qquad \mbox{ where }  \big|\{a,b,c,d\}\big| = 4 \mbox{ and } a<b \mbox{ and } c<d,\\
(15) &&
\{\z,\e_a,\d_{bc},\d_{de}\}, & 
	\qquad \mbox{ where } \big|\{a,b,c,d,e\}\big| = 5 \mbox{ and } b<c \mbox{ and } d<e.
}
These 135, found by brute force, form bases 
of the tangent spaces to the 45 tori figuring in \Cref{thm:XLV},
and each torus is reflected by the product of the four corresponding root reflections.

For example, the span $\R^4 \x \{\mspace{-1mu}0\}^{\mn 4}$ of 
$\{\g_{12},\d_{12},\g_{34},\d_{34}\}$
meets $\Phi^+$ in $\{\g_{ab}, \d_{ab} \colon 1 \leq a < b \leq 4\}$. 
Among these, the roots orthogonal to $\d_{ab}$
are $\{\g_{ab},\g_{cd},\d_{cd}\}$ (where $\big|\{a,b,c,d\}\big| = 4$)
and likewise the roots orthogonal to $\g_{ab}$ are $\{\d_{ab},\g_{cd},\d_{cd}\}$,
so the spanning quadruples are determined by the (three) partitions of $\{1,2,3,4\}$ into pairs of pairs $\big\{\{a,b\},\{c,d\}\big\}$.

\end{remark}

\appendix

\section{Leray and Koszul's theorem on $H^*(G/S^1)$}%
\label{sec:nosurj}

In order to obtain \Cref{theorem:main}, we needed some grasp on the 
cohomology ring $\H(G/S)$ of a homogeneous space $G/S$, 
for $G$ compact connected and $S$ a circle.

\thmHGS*

We found belatedly that this is a trivial generalization of long-known results.
General statements on the cohomology of 
a homogeneous space were already available to Leray in 1946, 
the year after his release from prison~\cite[\SS{3}, item~(4)]{miller2000leray}. 
In the second of his four \emph{Comptes Rendus} announcements from that year%
~\cite[bottom of p. 1421]{lerayCR1946b},
he states the following result.\footnote{ %
	See also Borel~\cite[par.~12]{borel1leraybio};
	only due to Borel's account are we confident 
	``compact Lie group'' was the accurate contemporary reading of Leray's \emph{groupe bicompact}.
}

\begin{restatable}[Leray, 1946]{theorem}{Leray}\label{Leray}
	Let $K$ be a compact, simply-connected Lie group
	and $S$ a closed, one-parameter subgroup [\emph{viz.} a circle].
	Write $\pi\: K \lt K/S$ for the projection.
	Then $\H(K/S;\Q)$ is generated as a commutative graded algebra 
	by finitely many classes $z_\a$ of odd degree
	and one class $s \in H^2(K/S;\Q)$, subject to the sole relation $s^{n+1} = 0$
	for a certain [positive natural] $n$.
	The ring $\H(K;\Q)$ is freely generated as a commutative graded algebra
	by the classes $\pi^* z_\a$ and one further class $z_{2n+1} \in H^{2n+1}(K;\Q)$.
\end{restatable}
More explicitly, 
if $P$ is a homogeneous vector space of generators
for the exterior algebra $\H K = \ext P$, 
then the image of $\H(K/S) \lt \H K$
is an exterior subalgebra $\ewP$ on a subspace $\wP \iso P/\Q z_{2n+1}$ 
of codimension $1$, and lifting $\wP$ back to $\H(K/S)$ 
induces a $\Q$-algebra isomorphism
\quation{\label{eq:Leray}
	\H(K/S) \,\iso\, \quotientmed{\Q[s]\,}{(s^{n+1})} \,\ox\, \ewP.
}
The second clause of \Cref{thmHGS} is clearly a refinement of this result;
if one omits Leray's hypothesis $K$ be simply-connected and 
admits the possibility $n$ be $0$,
then so is the first clause.

The following year, 
Koszul published a note~\cite[p.~478, display]{koszul1947homologie},
also in the \emph{Comptes Rendus}, 
regarding Poincar{\'e} polynomials for these spaces,
which implies $n=1$ in Leray's result.

\begin{restatable}[Koszul, 1947]{theorem}{Koszul}\label{Koszul}
	Let $K$ be a compact, connected Lie group and $S$ a compact, connected 
	$1$-dimensional subgroup [again, a circle]
	such that the image of $H_1(S;\Q) \lt H_1(K;\Q)$ is zero.
	Then the Poincar{\'e} polynomials (in the indeterminate $t$) 
	of $K/S$ and $K$ are related by
	\[
	p(K) (1 + t^2) = p(K/S)(1+t^3).
	\]
\end{restatable}
Koszul, unlike Leray, does include an indication of a proof,
which we translate without elaboration, 
leaving it to the reader to decide for themself
how much further detail they require
and provide it if they can.
After we will provide an alternate proof of \Cref{thmHGS}
and hence of Leray's and Koszul's theorems.

\begin{proof}[Koszul's proof]
	A choice of $K$-biinvariant Riemannian metric $B$ on $K$ 
	induces an isomorphism $\phi\:v \lmt B(v,-)$ from 
	the Lie algebra $\fk$, 
	conceived as the space of left-invariant vector fields on $K$,
	to the space $\Omega^1 \ceq \Omega^1(K)^K$ 
	of left-invariant $1$-forms.
	This allows us to define a Lie bracket on $\Omega^1(K)^K$,
	and to associate to $S$ the Lie subalgebra $\Omega^{1,0} \ceq \phi(\fs)$
	and its $B$-orthogonal complement $\Omega^{0,1}$.
	Then the differential algebra $\Omega^{\bl}(K)^K$ is bigraded by 
	$\Omega^{p,q} = \bigwedge^p \Omega^{1,0} \wedge \bigwedge^q \Omega^{0,1}$,
	and particularly we may consider the spectral sequence
	associated to the filtration by ideals $I^q = \Omega^{\bl,\geq q}$. 
	In this spectral sequence, one has
 	\eqn{
  		E_1^{0,\bul} &\iso \Omega^\bul(K/S)^K,\\
  		E_2^{0,\bul} &\iso \H(K/S),\\
 		\Ei \col &\iso \im\mn\big(\mn H^*(K/S) \lt \H K\big).
 	}
	Observe that given any nonzero element $\a \in \Omega^{1,0}$,
	we always have $d\a \in I^2$.
	We can uniquely decompose the Cartan invariant $3$-form
	$\w\: u \wedge v \wedge w \lmt B\big([u,v],w\big)$ on $K$ as
	$\w =\sum \w_j$
	for $\w_j \in \Omega^{3-j,j}$.
	Now $d\w = 0$ and $\w_0 = 0$, so we have
	\[
	(d\a)^2 = 
	d(\a \wedge d\a) = 
	3B(\a,\a) d\w_2 = 
	-3B(\a,\a) d\w_3,
	\]
	which simultaneously lies in $I^4$ and is the exterior derivative of 
	an element of $I^1$.
	Thus the image of $\H(K/S) \lt \H K$ 
	cannot contain the class $[\w]$.\footnote{\ 
		This is not made explicit by Koszul, 
		but we have $\Omega^{\geq 2,\bl} = 0$,
		so $\w = \w_2 + \w_3$ really.
		If we pick a $B$-orthonormal basis of $\Omega^1$ 
		including $\a$, 
		and expand in terms of structure constants, 
		we~\cite{MO:Bryant} get
		$\w_2 = \a \wedge d\a$.
	}
\end{proof}

Before our proof, 
we illustrate with a representative example 
the features of the general case.

\bex
Let $S$ be a circle contained in the second factor 
of the group $G = \U(2) \x \Sp(1)$. 
The cohomology of $G$ is the exterior algebra 
$\H G = \Lambda[z_1,z_3,q_3]$, where $\deg z_1 = 1$ and $\deg z_3 = \deg q_3 = 3$, 
and the cohomology $\HS = \H(BS) = \Q[s]$, where $\deg s = 2$.
Since $G/S = \U(2) \x \big(\Sp(1)/S\big) \homeo \U(2) \x S^2$,
we expect to find $E_\infty \iso \big(\Q[s]/(s^2)\big) \ox \Lambda[z_1,z_3]$
in the \SSS $(E_r,d_r)$ associated to $G \to {}_S G \to BS$.
Indeed, its $E_2$ page is the tensor product $\HS \ox \H G$. 
From the fact the map $H^1 G \longto H^1 S$ is zero
it will be shown to follow 
that the differential $d_2$ is zero. 
Next, $E_4 = E_2$ for lacunary reasons.
The differential $d_4$ can be shown to annihilate each of $s$, $z_1$, $z_3$ 
and take $q_3 \longmapsto s^2$. 
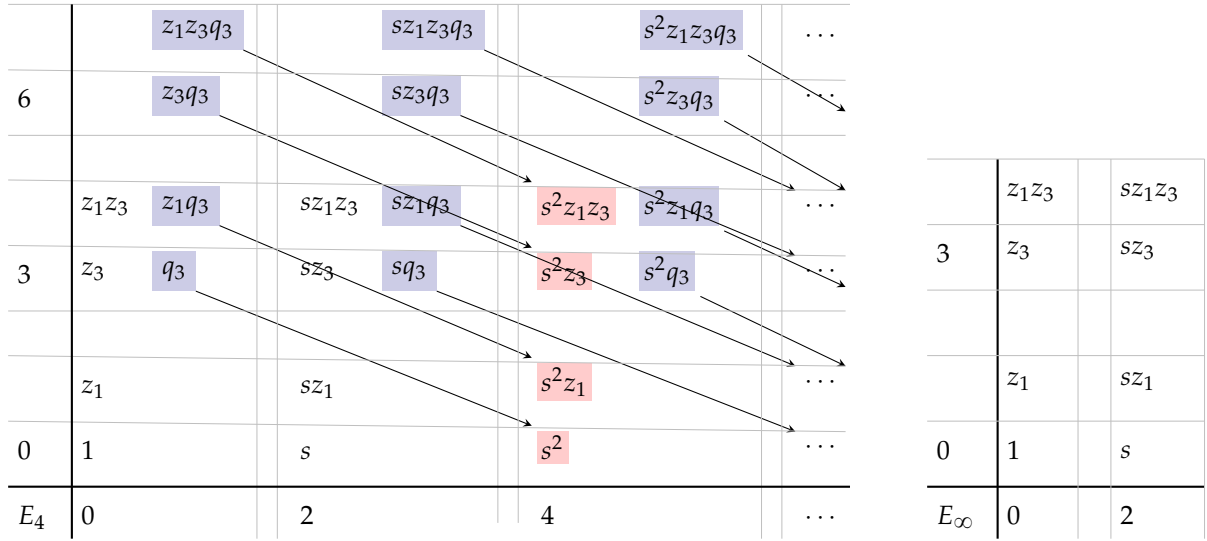
\[
\begin{tikzpicture}
\tikzstyle{every node}=[font=\small]
\matrix (m) [
matrix of math nodes,
nodes in empty cells,
nodes={
	anchor=west,
	minimum width=1.5ex,
	minimum height=1.5ex,
	outer sep=0pt
},
column sep=1ex,
row sep=1ex,
column 1/.style={anchor=base east},
column 2/.style={anchor=base east},
column 3/.style={anchor=base east},
column 4/.style={anchor=base east},
column 5/.style={anchor=base east},
column 6/.style={anchor=base east}
column 7/.style={anchor=base east}
column 8/.style={anchor=base east}
]
{
	\vpz		&\vpz& z_1z_3q_3	&\vpz	&& sz_1z_3q_3	&\vpz&& s^2 z_1z_3q_3&\vpz 	&\cdots \vpz\\
	\vpz	6	&& z_3q_3	&	&& sz_3q_3	&	&& s^2 z_3q_3	& 	&\cdots\\
	\phantom{\cdots}&&&&&&&	&&&\phantom{\cdots}\\
	\vpz		& z_1z_3&z_1q_3 	&	& sz_1z_3&sz_1q_3	&	& s^2z_1z_3&s^2z_1q_3	&	&\cdots\\
	\vpz	3	& z_3 &q_3&	& sz_3&sq_3	&	& s^2z_3&s^2 q_3	&	& \cdots\\
	\phantom{\cdots}&&&&&&&&&&\phantom{\cdots}\\
	\vpz		&z_1	\vpz	&&	& sz_1\vpz	&&	& s^2 z_1	&&	& \cdots\\
	\vpz	0	& 1	\vpz	&&	& s	\vpz	&& 	& s^2 	&&	&\cdots \\
	\vpz E_4	& 0		\vpz&&\vpz	&2	\vpz	&\vpz &	&\vpz 4 	 	\vpz&	\vpz&\vpz&\cdots\vpz\\
};
\path[-stealth] 	(m-1-3.south east) edge (m-4-8.north west)
(m-2-3.south east) edge (m-5-8.north west)
(m-4-3.south east) edge (m-7-8.north west)
(m-5-3.south east) edge (m-8-8.north west)
(m-1-6.south east) edge (m-4-11.north west)
(m-2-6.south east) edge (m-5-11.north west)
(m-4-6.south east) edge (m-7-11.north west)
(m-5-6.south east) edge (m-8-11.north west)		
(m-1-9.south east) edge (m-2-11.south east)
(m-2-9.south east) edge (m-4-11.north east)
(m-4-9.south east) edge (m-5-11.south east)
(m-5-9.south east) edge (m-7-11.north east);
\begin{pgfonlayer}{background}
\node [nonzero, fit=(m-1-3)] {};
\node [nonzero, fit=(m-2-3)] {};
\node [nonzero, fit=(m-4-3)] {};
\node [nonzero, fit=(m-5-3)] {};
\node [nonzero, fit=(m-1-6)] {};
\node [nonzero, fit=(m-2-6)] {};
\node [nonzero, fit=(m-4-6)] {};
\node [nonzero, fit=(m-5-6)] {};
\node [nonzerosm, fit=(m-1-9)] {};
\node [nonzerosm, fit=(m-2-9)] {};
\node [nonzerosm, fit=(m-4-9)] {};
\node [nonzerosm, fit=(m-5-9)] {};        
\node [image, fit=(m-4-8)] {};
\node [image, fit=(m-5-8)] {};
\node [image, fit=(m-7-8)] {};
\node [image, fit=(m-8-8)] {};
\end{pgfonlayer}
\path[lightgray, ultra thin]       
(m-9-4.south west) edge (m-1-4.north west)
(m-9-4.south east) edge (m-1-4.north east)
(m-9-7.south west) edge (m-1-7.north west)
(m-9-7.south east) edge (m-1-7.north east)
(m-9-10.south west) edge (m-1-10.north west)
(m-9-10.south east) edge (m-1-10.north east)
(m-8-1.north west) edge (m-8-11.north east)
(m-7-1.north west) edge (m-7-11.north east)
(m-6-1.north west) edge (m-6-11.north east)
(m-5-1.north west) edge (m-5-11.north east)
(m-4-1.north west) edge (m-4-11.north east)
(m-3-1.north west) edge (m-3-11.north east)
(m-2-1.north west) edge (m-2-11.north east)
(m-1-1.north west) edge (m-1-11.north east);
\path[thick] 	(m-9-2.south west) edge (m-1-2.north west)
(m-9-1.north west) edge (m-9-11.north east);
\end{tikzpicture}
\qquad
\begin{tikzpicture}
\tikzstyle{every node}=[font=\small]
\matrix (m) [
matrix of math nodes,
nodes in empty cells,
nodes={
	anchor=west,
	minimum width=.5ex,
	minimum height=1.5ex,
	outer sep=0pt
},
column sep=1ex,
row sep=1ex,
column 1/.style={anchor=base east},
column 2/.style={anchor=base east},
column 3/.style={anchor=base east},
column 4/.style={anchor=base east},
]
{
	\vpz		& z_1z_3	\vpz	&	\vpz	& sz_1z_3			\vpz	&\vpz\\
	\vpz	3	& z_3		&		& sz_3\phantom{z_1}	&\vpz\\
	\vpz		&			&		&\phantom{sz_1z_3}		&\vpz\\
	\vpz		& z_1		&		& sz_1\phantom{z_3}	&\vpz\\
	\vpz	0	& 1\vpz		&		& s\phantom{z_1z_3}\vpz	&\vpz\\
	\vpz E_\infty	& 0	\vpz	&	\vpz	&2\phantom{z_1z_3}	\vpz	&\vpz\\
};
\path[thick] 	(m-6-2.south west) edge (m-1-2.north west)
(m-6-1.north west) edge (m-6-5.north west);
\path[lightgray, ultra thin]       
(m-6-3.south west) edge (m-1-3.north west)
(m-6-4.south west) edge (m-1-4.north west)
(m-6-5.south west) edge (m-1-5.north west)
(m-5-1.north west) edge (m-5-5.north west)
(m-4-1.north west) edge (m-4-5.north west)
(m-3-1.north west) edge (m-3-5.north west)
(m-2-1.north west) edge (m-2-5.north west)                      
(m-1-1.north west) edge (m-1-5.north west);
\end{tikzpicture}
\]
Because $d_4$ is an antiderivation, its kernel is the subalgebra $\Q[s] \ox \Lambda[z_1,z_3]$
and its image the ideal $(s^2)$ in that subalgebra.
Elements mapped to a nonzero element by $d_4$ are marked as blue in the diagram
and elements in the image in red; the vector space spanned by these elements vanishes in $E_5$.
Thus $E_5 = \big(\Q[s]/(s^2)\big) \ox \Lambda[z_1,z_3]$.
For lacunary reasons, $E_5 = E_\infty$.
\eex

We work with a general compact, connected Lie group $G$
and closed, connected subgroup $H$,
specializing to the desired case at the end.
Because the Borel fibration $G \to G_H \to BH$
is a principal $G$-bundle, 
it admits a classifying map to $BG$,
which can be seen to be (homotopic to) 
the map $BH = EG/H \to EG/G = BG$
functorially induced by the inclusion $H \inc G$.
The resulting map of principal $G$-bundles
\[
\xymatrix@C=1.75em{
	G  \ar@{=}[r]\ar[d]	&G\ar[d]\\
	G_H\ \ar[r]_\psi	\ar[d]	&EG \ar[d]\\
	BH \ar[r]_{\rho}			&BG
}
\]
induces a map $(\psi^*_r)$ of Serre spectral sequences.
Each page of the right sequence $(\tE_r,\td_r)$ is of tensor form,
and the transgressions $\td_{2k}\: \tE_{2k}^{0,2k-1} \lt \tE_{2k}^{2k,0}$
induce~%
\cite[Thm.~13.1]{borelthesis}
a degree-one linear isomorphism
\[
P\H G \,\isoto\, \quotientmed{H_G^{\geq 1}}{H_G^{\geq 1}H_G^{\geq 1}}
\]
between the space of primitive elements of the Hopf algebra $\H G$
and the space of indecomposables of the polynomial ring $\HG = \H BG$,
which one should think of as residues of homogeneous generators.
This bijection completely determines the differentials $\td_r$,
and in turn the differentials of the left spectral sequence $(E_r,d_r)$
are completely determined by the chain relations $\psi^*_r \td_r = d_r \psi^*_r$.
A lifting of the linear isomorphism to a degree-one linear injection 
$\tau\: P\H G \lt \H BG$, followed by the map $\rho^*\: \H BG \lt \H BH$
induces a unique derivation $d = \rho^* \o \tau$ 
on the page $E_2 = \H BH \ox \H G$
which vanishes on the $\H BH$ factor 
and simultaneously lifts all the differentials $d_r$.
Borel shows~\cite[Thm.~25.2]{borelthesis} 
the cohomology of the resulting algebra 
$(\H_H \ox \H G,d)$,
the \defd{Cartan algebra} 
is isomorphic to $\H(G/H)$.\footnote{\ 
	Cartan earlier arrived at the same algebra by very different methods%
	~\cite[Thm.~5, p.~216]{cartan1950transgression}.
	Borel's proof can be seen in retrospect to be a consequence of a general
	method in rational homotopy theory~\cite[Prop.~15.5,8]{FHT}
	which converts 
	compatible models
	of a fibration $E \lt B$ and of a map $\rho\: B' \lt B$
	into a model of the total space of the pullback $\rho^* E \lt B'$.
}
In fact, one recovers the \SSS again 
as the spectral sequence
induced from the Cartan algebra
by the filtration induced from the grading of $\H_H$.

An important feature of this \CDGA
is that typically some of the differentials $dz$ of primitives $z \in P\H G$
are ``redundant'' in the sense they lie in the ideal $H_G^{\geq 1} 
\. d(P \H G)$ 
generated by positive-degree multiples of such differentials. 
The space $\defm{\wP}$ of these primitives 
with redundant differential is called the \defd{Samelson space},
and if we denote its complement by $\defm{\cP} \ceq P \H G / \wP$,
the filtration spectral sequence induced by the grading on $\HK$ 
shows the Cartan algebra factors as a tensor product
\[
(\H_H \ox \H G,d) \iso (\H_H \ox \ecP,d) \ox (\ewP,0)
\]
of {\CDGA}s {\cite{andre1962tohoku}%
	\cite[Thm.~2.15.V, p.~73]{GHVIII}%
	\cite[Prop.~8.5.4, p.~141]{onishchik}}; 
moreover, 
viewing the filtration spectral sequence as the \SSS of $G \to G_H \to BH$,
we may identify $\ewP$ with the image of $\H(G/H) \lt \H G$.
A \defd{pure Sullivan algebra} 
is a free commutative graded algebra $\Q[Q] \ox \ext P$
on an evenly- and positively-graded rational vector space $Q$
and an oddly- and positively-graded $P$
equipped with a derivation $d$ vanishing on $Q$ such that $d^2 = 0$
and $dP < \Q[Q]$.
For such a \CDGA, a Samelson space $\wP$ is similarly defined 
as $\{z \in P : dz \in (Q \.dP)\}$.

\bpf[Proof of \Cref{thmHGS}]
If $H^1 G \lt H^1 S$ is surjective, 
then $\H G \iso \H S \ox \H(G/S)$ 
by Samelson's theorem~\cite[Satz VI(b), p. 1134]{samelson1941samelson},
yielding the first clause. 

Otherwise, we compute $\H(G/S)$
as the cohomology of the Cartan algebra $(\HS \ox \H G,d)$.
Write $\HS = \Q[s]$ for $s \in H^2 BS^1$.
Since $\Q[s]$ is a graded principal ideal domain,
in any homogeneous basis $(z_j)$ of $PG$, all but one $dz_j$
is a redundant generator of the ideal $(dz_j) \ideal \Q[s]$,
so the Samelson subspace $\smash\wP$ generating $\im\mn\big(\mn\H(G/S) \to \H G\big)$ 
has dimension $\rk G - 1$,
and hence $\H(G/S)$ has the form claimed in (\ref{eq:Leray})
(i.e.,
$G/S$ is formal 
in the sense of rational homotopy theory).
The map $H^2 BG \lt H^2 BS = \Q \. s$ 
is conjugate through transgression isomorphisms to the map $H^1 G \lt H^1 S$
and hence by assumption is trivial.
It follows from \Cref{theorem:Hsurjecttorus} 
that $S$ lies in the commutator subgroup
$K$ of $G$
and we can factor the map of interest as
${\HG \to \HK \to \HS}$.
The first map is surjective since $G$ has a finite central extension $\wt G$
of the form ${\wt K \times (\wt G/\wt K)}$,
so that $\smash{\H_{\wt K} \iso \HK}$ 
is a tensor factor of $\smash{\H_{\wt G} \iso \HG}$
(\Cref{pathtriv})
and we may just consider the image of $\HK \lt \HS$.
By the following lemma 
this is $(s^2)$,
so $n = 1$ in (\ref{eq:Leray}).
\epf

\begin{lemma}\label{imageHBK}
	Let $K$ be a semisimple Lie group containing a circle $S$.
	The image of $\HK \longto \HS \iso \Q[s]$ contains $s^2 \in H^4_S$.
\end{lemma}
\begin{proof}
	Let $T$ be a maximal torus of $K$ containing $S$. 
	By \Cref{KSdecomp},
	$\HK \lt \HT$ is an injection with image the invariant subring $(\H_T)^W$ 
	under the action of $W = W_K$.
	Write $\R[\f t]$ for the graded algebra of polynomial functions 
	on the Lie algebra $\ft$ of $T$, 
	assigning nonzero linear forms degree $2$.
	Extending coefficients to $\R$,
	the Chern--Weil homomorphism~\cite[Thm. 2.4]{kobayashinomizuII}
	and Chevalley restriction theorem \cite[\SS{{IV}}]{chevalley1950betti}
	translate the sequence 
	\[
	\HK \isoto (\HT)^W \longinc \HT \longepi \HS
	\]
	into 
	\[
	\R[\fk]^K \xrightarrow[\substack{\sim\\ \vphantom{X}}]{\textrm{rest}} 
	\R[\ft]^W \longinc 
	\R[\ft] \xepi{\textrm{rest}} 
	\R[\fs].
	\]
	In particular, elements of $H^4_K$
	correspond to $W$-invariant quadratic forms on $\f t$
	and $H^4_K \lt H^4_S$ is surjective if
	any such form does not vanish on $\f s$.
	But the Killing form $B$ of $K$
	is a $(\Ad K)$-invariant bilinear form on $\fk$,
	negative definite
	since $K$ is semisimple~\cite[Prop.~V.(5.13), p.~214]{brockertomdieck},
	so precomposing the diagonal inclusion 
	$\smash{
		\f t \mono \f t^2 \inc \f k^2
	}$ 
	yields a $W$-invariant quadratic form on $\f t$
	restricting nontrivially to any one-dimensional subspace $\f s$.
\end{proof}

\brmk
The author's original proof of this lemma proceeded laboriously 
by cases through the simple groups. He is indebted to Mathew Wolak
for pointing out the Killing form is invariant and definite.
\ermk

\section{Partial reductions}\label{sec:failures}

Some fragments of the results we are interested in
persist even in the case $G$ is merely assumed 
to be a pro-Lie group, not necessarily connected,
but as the surviving results are not so powerful 
as one might like, they have been deferred to this appendix. 
We can nevertheless prove the expected result when the isotropy 
group remains a circle.

\subsection{Connected groups}\label{sec:connected}
Let $G$ be a topological group and $K$ a closed subgroup. 
We would like to reduce the question of when $(G,K)$ 
is isotropy-formal to 
the same for connected components $(G_0, K_0)$
of the identity in each,
but that is too much to hope. 
There is at least the following diagram:
\[
\xymatrix@C=1.25em@R=3.25em{
\ \ \, G_0/K_0   \ar@{}[r]^(.325){}="a"^(.85){}="b"\ar"a";"b" 		  \ar[d]^i	&
\ \ \, G/K_0	 \ar@{}[r]^(.325){}="a"^(.85){}="b"\ar"a";"b"^(.35)\d \ar[d]^j	&
\ G/K  				\ar[d]^k	&\\
 {}_{K_0}\mspace{-1.5mu} G_0/K_0  	\ar[r]	  & 
 {}_{K_0}\mspace{-1.5mu} G/K_0  	\ar[r]^\e \ar[d]^\h& 				
 {}_{K_0}\mspace{-1.5mu} G/K 		\ar[d]^\t 	\\ &
\ \ {}_K G/K_0		\ar@{}[r]^(.325){}="a"^(.8){}="b"\ar"a";"b"^(.4)\z&
\ \,{}_K G/K.
}
 \]
As $K_0$ lies in $G_0$, the map $j$
can be understand as the disjoint union of $\pi_0 G$ parallel copies of $i$,
so the one surjects in cohomology just if the other does.
Less can be said about the other maps.

\bprop\label{thm:conn}
Assume $\pi_0 K $ is finite. 
If $j^*\: \H_{K_0}(G/K_0) \lt \H(G/K_0)$ is surjective, 
then so is $k^*\:\H_{K_0}(G/K) \lt \H(G/K)$,
and $(\t \o k)^*$
is surjective as well if and only if additionally the left action of $K$ on $G/K$ 
induces a trivial action of $\pi_0 K$ on $\H(G/K)$.
Suppose additionally $K$ lies in $G_0$ 
and $\HKz$ is free over $\HK$.
Then if $k^*$ is surjective, so also are $(\t \o k)^*$
and $(\h \o j)^*$ and $j^*$.
\eprop
\bpf
The maps $\d$, $\e$, and $\z$ in the diagram are
coverings induced by a right $\pi_0 K$-action
in such a way that $j$ and $\eta$ and hence $j^*$ 
and $(\eta \o j)^*$ are $\pi_0 K$-equivariant.
Since we assume $\pi_0 K$ is finite,
a standard lemma \cite[Prop.~3G.1]{hatcher} 
identifies $\d^*$, $\e^*$, and $\z^*$ with inclusions of invariants
so that $k^*$ becomes the restriction 
$\H_{K_0}(G/K_0)^{\pi_0 K} \lt \H(G/K_0)^{\pi_0 K}$
and $(\t \o k)^*$ the restriction
$\H_K(G/K_0)^{\pi_0 K} \lt \H(G/K_0)^{\pi_0 K}$.
If $j^*$ is surjective, 
then $k^*$ must be as well,
and if $(\h \o j)^*$ is, then so is $(\t \o k)^*$, 
in both cases by averaging. 
Now, the map $(\t \o k)^*\:\HK(G/K) \lt \H(G/K)$ 
is surjective if and only if the \SSS associated to the Borel fibration $G/K \to {}_K G/K \to BK$ collapses at $E_2$ 
and the action of $\pi_1 BK$ on the cohomology of the fiber $G/K$ is trivial~\cite[Prop.~4.1, p.~129]{borelthesis}
so triviality of the action is necessary.
On the other hand, if $k^*$ is surjective and the action is trivial,
then the map of {\SSS}s induced by the map ${}_{K_0} G/K \epi {}_{K} G/K$ 
is represented on the $E_2$ page 
by an injection $\HK \ox \H (G/K) \longmono \HKz \ox \H (G/K)$, 
and since all differentials vanish on 
$\smash{E_2\col \iso \Q \ox \H (G/K)}$ in the larger
sequence, the same holds in the smaller, 
so it also collapses and $(\t \o k)^*$ is surjective. 
%

In general, in a Borel fibration $X \to X_K \to BK$,
the action of $\pi_1 BK = \pi_0 K$ on the fiber $X$
descends from the action of $K$ on $X$,
so if we assume $K$ lies in $G_0$, 
then by path-connectedness of the latter, $\pi_0 K$
acts trivially on the right on the fibers 
$G$, 
${}_{K_0}\mspace{-1.5mu} G$, and
${}_K G$
of the Borel fibrations 
over $BK$,
the cohomology of whose total spaces is in question.
If we assume additionally that $\HKz$ is free over $\HK$,
then \Cref{freeECext}
applies to identify $j^*$ with $\smash{\id_{\HKz} \ox_{\HK} \,k^*}$
and $(\h \o j)^*$ with  $\smash{\id_{\HKz} \ox_{\HK} \,(\t \o  k)^*}$,
meaning in either pair of maps, the latter is surjective 
if and only if the former is.
If $k^*$ is surjective, then, 
by the argument of the previous paragraph, 
so also is $(\t \o k)^*$, and then by \Cref{freeECext}
so also are $(\h \o j)^*$ and $j^*$.
%
\epf

As limiting as the hypotheses seem, they are necessary.
We will discuss their disappointing asymmetry in
\Cref{connectlimit}.

\blem\label{freeSerreext}
Let a map of fibrations with homotopy fiber $F$
be given as in (\ref{overbundlemap})
such that $\pi_1 B_0$ acts trivially on $H^* F$
and $H^* B$ is a flat module over $H^* B_0$.
Then there is an $\H E_0$-algebra isomorphism
\[
	\smash{\defm{\psi}\: \xt{\H B_0}{\H B\!}{\!\H E_0} \isoto \H E}
\] 
natural in $\xi$. 
\elem
\bpf
The map induces a map $(\defm{\psi_r^0})$ 
of Serre spectral sequences $(E^0_r,d^0_r) \lt (E_r,d_r)$.
As each $E^0_r$ is an $\H B_0$-algebra and each $E_r$ an $\H B$-algebra,
we obtain a collection of maps 
\[
	\defm{\psi_r}					\:
	\defm{E'_r} 						\ceq
	\xt{H^* B_0}{\H B\!}{\!E^0_r} 	\lt 
	\xt{H^* B}{\H B\!}{\!E_r} 		\isoto
	E_r.
\] 
If we assign $E'_r$ the differential $\defm{d'_r} \ceq \id \ox\, d^0_r$,
then $(E'_r,d'_r)$ is a spectral sequence by flatness:  
\[	
	\H E'_r = 
	\xt{\H B_0}{\H B\!}{\!\H E_r^0} = 
	\xt{\H B_0}{\H B\!}{\!E_{r+1}^0} = 
	E'_{r+1}.
\]
Since $(\psi_r^0)$ was a spectral sequence map, 
so also is $(\psi_r)$.
As we assume simple coefficients, $\psi_2$ is the canonical isomorphism.
Inductively, since each cochain map $\psi_r$ is an isomorphism, 
so also is the map $\psi_{r+1}$ it induces in cohomology.
Thus $\psi_\infty$ is an isomorphism.
As $\psi_\infty$ is the map of associated graded algebras 
induced from $\psi$,
it follows $\psi$ is an isomorphism as well.
\epf

\bcor\label{freeECext}
Let a Lie group $K$ act on $X$
in such a way that the action of $\pi_1 BK$ on $\H X$
induced by the Borel fibration $X \to X_K \to BK$ is trivial.
Suppose $H$ is a subgroup of $K$ such that
$\H_H$ is free as an $\HK$-module.
Then there is an isomorphism $\H_H \ox_{\HK}\HK X \isoto \H_H X$
natural in $X$.
\ecor
\bpf
Apply \Cref{freeSerreext} to the map $X_H \lt X_K$.
\epf
%

\brmk\label{connectlimit}
In case $\HK$ is not free over $\HG$,
\Cref{freeECext} can fail. 
To see this,
consider the block-diagonal inclusion of
$H = \SU(3)^2$ in $K = \SU(6)$ and let each act on the right of $X = \U(6)$
by multiplication.
We want to determine whether the map
\[
\xt{\H_{\SU(3)^2}}{\H_{\SU(6)}}{\H\big(\mspace{-1mu}\U(6)/\SU(6)\big)}
	\lt
\H\big(\mspace{-1mu}\U(6)/\SU(3)^2\big)
\]
is an isomorphism.
But $\U(6) /\SU(6) \iso S^1$
has cohomology ring $\ext z_1$ concentrated in odd degree,
so the map
$\H_{\SU(6)} \lt \H\big(\mspace{-1mu}\U(6) / \SU(6)\big)$
is trivial
and the domain is isomorphic to 
\[
	\big(\mn\xt{\H_{\SU(6)}}{H_{\SU(3)^2}}{\Q} \big)\ox \Lambda z_1.
\]
On the other hand, 
it is easy to see from the Cartan algebra
of \Cref{sec:nosurj} that
$\H\big(\mspace{-1mu}\U(6) / \SU(3)^2\big) \iso
\H\big(\SU(6) / \SU(3)^2\big) \ox \Lambda z_1$,
so the map in question is an isomorphism only if 
${H_{\SU(3)^2} \ox_{\H_{\SU(6)}} \Q}
\lt \smash{\H\big(\SU(6) / \SU(3)^2\big)}$
is.
This is, however, untrue \cite[pp.~486--488]{GHVIII}:
the target is the ring 
${\Tor^*_{\H_{\SU(6)}}\mn(\Q,\H_{\SU(3)^2})}$
and the sources the proper subring 
\smash{${\Tor^0_{\H_{\SU(6)}}\mn(\Q,\H_{\SU(3)^2})}$}.

The condition that $K$ lie within $G_0$ is severe as well,
but without it, the right action of $K$ on $G/K_0$
already induces a nontrivial action of $\pi_1 BK = \pi_0 K$
on $H^0(G/K_0)$.
\ermk

\subsection{Lie groups}\label{sec:Lie}

To make as complete as possible 
the attempted reduction of the problem of isotropy-formality
to the case of a torus in a semisimple group,
we include the case of compact groups.
We get surprisingly far, as there are relatively
few algebraic obstacles,
but we only achieve a complete reduction if 
the isotropy group is Lie.
In case the isotropy group is a circle,
we do get back a version of \Cref{algor}, 
namely \Cref{thm:cptoverlie}.

Every compact Hausdorff group $G$ can be realized as an 
{inverse limit} of Lie group homomorphisms%
~\cite[Ex.~3.4, p.~137]{hofmannmorrispro}, 
which is to say the limit in the category of topological groups
of a directed system 
\[
	(G_\a,\,\phi_{\a \to \b}\: G_\a \lt G_\b)_{\a \geq \b}
\] 
of Lie groups, 
the maps $\phi_{\a \to \b}$ between which may be taken surjective~%
\cite[Prop.~1.33, p.~21]{hofmannmorriscompact}.
Such a realization comes equipped with unique surjections 
$\defm{\phi_\a}\: G \lt G_\a$ for each $G_\a$ 
such that $\phi_\b = \phi_{\a \to \b} \o \phi_\a$ whenever $\a \geq \b$.
If $K$ is a closed subgroup of $G$, 
let $\defm{K_\a} \ceq \phi_\a K \leq G_\a$;
then the restrictions $\rest{\phi_{\a \to \b}}{K_\a}$ 
realize $K$ as $\limit K_\a$. 
The inclusion map of inverse systems 
$(K_\a,\rest{\phi_{\a\to\b}}{K_\a}) \lt (G_\a,\phi_{\a \to \b})$ 
induces a quotient system $(G_\a/K_\a,\bar \phi_{\a \to \b})$ 
of continuous surjections of homogeneous spaces
and the left action of $(K_\a,\rest{\phi_{\a\to\b}}{K_\a})$
induces a system $({}_{K_\a}\mn G_\a/K_\a,\bar\phi'_{\a\to\b})$ 
of homotopy quotients.
The canonical map 
$G/K \lt \lim G_\a/K_\a$ is a continuous bijection of compact Hausdorff spaces,
hence a homeomorphism 
(in fact, this is still the case if $G$ and $K$
are non-compact pro-Lie groups~\cite[Lem.~1]{mostert1953local}).
We take as our realization of $E(-) \lt B(-)$ 
the Milnor construction~\cite{milnor1956universalII}.
The functorially induced map $EG \lt \lim EG_\a$
is actually a $G$-equivariant homeomorphism,
inducing a homeomorphism $BG \lt \lim BG_\a$.
Thus the map $EK \x G/K \lt \lim (EK_\a \x G_\a/K_\a)$
is a $K$-equivariant homeomorphism as well, so finally 
we can write ${}_K G/K$ as $\lim {}_{K_\a} \mn G_\a/K_\a$. 
Then the fiber inclusion $i\: G/K \lt {}_K G/K$ 
is identified with $\lim (i_\a\:G_\a/K_\a \lt {}_{K_\a}\mn G_\a/K_\a)$.

\v{C}ech cohomology 
(with coefficients in the constant sheaf $\ul\Q$, henceforth) 
converts inverse limits to direct limits \cite[pp.~318--9]{spanier}; 
the essential point is that an inverse 
limit can be viewed as an intersection. 

\bex\label{eg:sol}
The solenoid $\Xi$ which is the inverse limit of the sequence 
$\cdots \to S^1 \os 2 \to \smash{S^1 \os 2 \to S^1}$, 
though connected, 
has continuum-many path components,
so particularly $H^0 \Xi$ is large.
Nevertheless,
applying 
$\smash{\wH^0}$
to the sequence yields isomorphisms
$\smash{\cdots \from \Q \os \id\from \Q \os \id\from \Q}$
and 
$\smash{\wH^1}$ isomorphisms
$\smash{\cdots \from \Q \os 2\from \Q \os 2\from \Q}$,
so 
$\smash{\wH^0 \Xi \iso \Q \iso \wH^1 \Xi}$.
If we identify the map $\smash{\wH^1 S^1 \isoto \wH^1 \Xi}$
induced by projection to the last circle with $\smash{\id_\Q}$,
then projection to the $\smash{n\th}$-from-last induces multiplication by $1/2^n$.
\eex

Thus we can identify the restriction 
$\Cech(G/K \longinc {}_K G/K)$ 
 with 
%
%
\quation{\label{eq:colimCech}
\smash{
	\colim \mn\big(\mspace{-1mu}H^*_{K_\a} \mn(G_\a/K_\a) \xrightarrow{i^*_\a} 
	H^* (G_\a/K_\a)\big).
}
}
The following is then clear.

\bprop\label{thm:colimboring}
If there is a cofinal subset of indices $\a$ 
such that the associated $i_\a^*$ are surjective, then so is $i^*$. 
\eprop

But $i^*$ can be surjective though no individual 
map $H^*_{K_\a} \mn(G_\a/K_\a) \lt H^* (G_\a/K_\a)$ be.

\bex\label{EG:badEG}
Set $H_1 = \SU(6)$ and 
for each $k \geq 2$ set $H_k = \mathrm S\big(\U(3) \x \U(6)\big)$.
Let $G$ be the product $\prod_{k \geq 1} H_k$
and $K$ the subgroup 
$\big\{(A_1\+B_1) \concat (B_{k-1},A_k\+B_k)_{k \geq 2} \in H_1 \x \prod_{k \geq 2} 
H_k : A_k, B_k \in \SU(3)\big\}$,
where $A_k \+ B_k \in \SU(6)$ denotes the $6\x 6$ block-diagonal matrix with nonzero
$3\x 3$ blocks $A_k$ and $B_k$.
Then $(G,K)$ is isotropy-formal,
and is the limit of the quotients $G_n = 
\prod_{k \leq n} H_k$,
with the expected projections
$\phi_n \: G \lt G_n$ and $K_n = \phi_n K$,
but none of the pairs $(G_n,K_n)$ is isotropy-formal.

There is an evident artifice to this example.
The groups $H_k$ for $k \geq 2$
contain subgroups $H'_k = \SU(3) \x \SU(6)$
and also admit $\wt H_k = \SU(3) \x \SU(6) \x S^1$ as six-fold covers,
and these are decomposable.
Replacing $G$ with $\wt G = H_1 \x \prod_{k \geq 2} \wt H_k$,
with $\wt G_n = \prod_{k \leq n} \wt H_k$,
and $K$ with the isomorphic subgroup $\wt K$ of $G'$
with entries $1$ in all $S^1$ factors and $A_k,B_k$ 
in special unitary factors as before, 
or replacing $G$ with $G' = H_1 \x \prod_{k \geq 2} H'_k$ 
and maintaining the old $K$,
the cohomological behavior of 
$\big(\mn\wt G_n, \im(\wt K \to \wt G_n)\big)$ is the same as before,
each $\wt G_n \lt G_n$ being a $6^{n-1}$-fold central cover,
and the behavior of $\big(G'_n,\im (K \to G'_n)\big)$ 
is similar except that all the $H^* S^1$ 
tensor factors are lost.
But $\wt G$ is also the limit of the groups
$
	\wt G'_n = \big(\mspace{1mu}\SU(6) \x \SU(3) \x S^1\big)^n,
$
and $G'$ of the groups $G''_n = \big(\SU(6) \x \SU(3)\x\{1\}\big)^n$,
and the images $\wt K'_n$ of $\wt K \to \wt G'_n$
and $K'_n$ of $K \to G'_n$
are both isomorphic to 
$\big\{(A_k\+B_k,B_k,1)_{k \leq n} : A_k,B_k \in \SU(3)\big\}$,
so the pairs $(\wt G'_n,\wt K'_n)$
and $(G'_n,K'_n)$ \emph{are} all isotropy-formal.
Thus in a sense we only obtained this counterexample 
by perversely choosing a bad inverse system when better---%
up to finite coverings---were plainly available.
The author still does not know if more meaningful counterexamples exist.
\eex

In the event $G$ and $K$ are connected,
the pure Sullivan models of Cartan and Kapovitch~\cite[Thm.~5, p.~216]{cartan1950transgression}%
\cite[Thm.~25.2]{borelthesis}%
\cite[Prop.~1]{kapovitch2002biquotients}%
\cite[Thm.~3.50]{FOT}
express each $i_\a^*$ from (\ref{eq:colimCech})
as the map induced in cohomology by \CDGA maps
\quation{
	\label{eq:CartanKapovitch}
	(\H_{K_\a} \ox \H_{K_\a} \ox \H G_\a, \td_\a\big) \lt (\H_{K_\a} \ox \H G_\a,d_\a\big).
}
With some care, we can realize $i^* = \colim i_\a^*$
as the cohomology of a colimit of these models.

\bprop\label{thm:dumbmodel}
	Let $(G,K)$ be a pair of compact, connected Hausdorff groups.
	Then the cohomology of
	the fiber inclusion $G/K \lt \KGK$ 
	is induced by a map
	\quation{\label{eq:Cechcolimspace}
		(\Cech_K \ox \Cech_K \ox \Cech G,\td\big)	
			\lt 
		(\Cech_K \ox \Cech G, d).
	}
	of pure Sullivan algebras given as follows.\footnote{\
		It is tempting to call these algebras Sullivan models,
		but to do so would require \CGA quasi-isomorphisms
		from our algebras to the algebras of polynomial differential forms
		$\APL({}_K G/K)$ and $\APL(G/K)$.
		We can construct such maps at each level 
		of the inverse system,
		but as $\APL$ computes \emph{singular} cohomology,
		it is unclear we will still have quasi-isomorphisms
		when we are done.
	}
	The differential $d$ is the unique derivation vanishing on $\HK$ 
	and extending the composition
	\[\smash{
		P\Cech G \xrightarrow[\sim]{\tau} 
		\defm{Q}\Cech_G \eqc {\wH_G^{\geq 1}} \,/\, 
		\wH_G^{\geq 1}\wH_G^{\geq 1}
		\os{s}\lt \Cech_G \xrightarrow{\rho^*} \Cech_K,
	}\]
	where $\tau$ is the transgression 
	in the \SSS of $G \to EG \to BG$, 
	the map $s$ is a certain graded linear lifting 
	of the indecomposables of $\Cech_G$ to generators,
	and $\rho = B(K \inc G)$ is the canonical map.
	The differential $\smash{\td}$ 
	is the unique derivation 
	vanishing on $\smash{\Cech_K \ox \Cech_K}$ and taking $z \in P\Cech G$
	to $1 \ox dz - dz \ox 1 \in \Cech_K \ox \Cech_K$.
	The map of differential graded algebras is that quotienting 
	out the ideal $\wH_K^{\geq 1} \ox \Cech_K \ox \Cech G$.
\eprop

\bpf
First we show this is a map of pure Sullivan algebras,
then that it computes the map in cohomology claimed. 
For the former, we need only see 
the commutative graded algebras underlying the proposed models are free,
which is to say $\Cech_K$ is a polynomial ring 
and $\Cech G$ an exterior algebra.
This results from the rather restricted nature of
surjective homomorphisms $G \epi G'$
between compact, connected Lie groups:
such a map induces a surjection
$\f g \epi \f g'$ of reductive Lie algebras, 
which is a factor projection. 
The group map is thus finitely covered by a factor projection:
\[
	\xymatrix@C=1em@R=2.25em{
		G'' \x \tG' \ar@{->>}[r]\ar[d] 	& \tG' \ar[d]\\
		G \ar@{->>}[r]						&  G'.
	}
\]
If $K$ is a subgroup of $G$ and $K'$ its image in $G'$,
then $K \epi K'$ is 
likewise finitely covered by a factor projection 
$K'' \x \tK' \epi \tK'$.\footnote{\ 
But not necessarily in such a way that $K''$ is contained in $G''$
and $\tK'$ in $\tG'$. 
For example, let $K = \Delta K'$ 
be a diagonally embedded copy of $K < G$ in $G \x G$
and consider the factor projection $G \x G \epi G$.
}
Since we take rational coefficients, by \Cref{universalH}
the maps in cohomology are then tensor factor inclusions of the form 
$\H G' \simto \H \tG' \lt \H G' \ox \H G'' \simto \H G$
and $\H BK' \lt \H BK' \ox \H BK'' \simto \H BK$.
Thus each of the maps 
between the models of ${}_{K_\a}\mn G_\a/K_\a$
and $G_\a/K_\a$ may be replaced with a tensor factor inclusion.
But the direct limit of such a system is a tensor product by definition,
and a tensor product of free commutative algebras is again free.
As \v{C}ech cohomology converts inverse limits to direct limits,
we can substitute $\Cech_K$ for $\colim \H_{K_\a}$
and $\Cech G$ for $\colim \H G_\a$.

To see the cohomology of the map is as claimed,
we must construct the differentials 
to be the colimit of the Cartan 
and Kapovitch differentials for the Lie pairs $(G_\a,K_\a)$.
Note that these models are not quite functorial,
in the sense that the differentials $d_\a$ and $\td_\a$---%
given as in the statement of the theorem
if $(G,K) = (G_\a,K_\a)$---%
each depend on an arbitarily chosen section 
$Q \H_{G_\a} \lt \H_{G_\a}$
of the reduction $\H_{G_\a} \lt Q\H_{G_\a}$.
For there to be a colimit at all,
the sections must be chosen coherent in the sense the obvious squares
\quation{\label{eq:sectionsquare}
	\begin{gathered}
\xymatrix@C=3em@R=3.75em{
	Q\H_{G_\a}  \ar@{ >->}[d]_{s_\a} 
	& Q\H_{G_\b} \ar@{ >->}[d]^{s_\b} \ar@{ >->}[l]_{Q(B\phi_{\a \to \b})^*}\\
	H_{G_\a}& \H_{G_\b} \ar@{ >->}[l]^{(B\phi_{\a \to \b})^*} 
}
	\end{gathered}
}
commute for all $\a \geq \b$.
One might hope to achieve this by defining $s$
first and then restricting, but 
then it is not necessarily the case
that the image of the composition 
$Q \H_{G_\a} \mono Q\Cech_G \os s\to \Cech_G$
lies in the image of $\smash{H_{G_\a} \longmono \Cech_G}$.
Instead, note that
the limit $G$ will not change if we extend the diagram to include 
all quotients of all $G_\a$, so we do.
Next, since a finite covering induces an isomorphism in rational cohomology, 
we may, by picking one ring in each 
isomorphism class, replace the diagram of graded rings $\H_{G_\a}$ 
by a skeleton 
in which no nonidentity arrow is an isomorphism.\footnote{\ 
	To appreciate how drastic this reduction is,
	note that if the solenoid $\Xi$ 
	of \Cref{eg:sol} is a quotient of $G$,
	say $G = H \x \Xi$, 
	then the corresponding parts of the diagram of $\Q$-algebras
	comprise solely factors $\H_{H_\a}$ and $\H_{H_\a} \ox \H_{S^1}$.
	Of course, this staggering swindle is only possible because we have 
	already	passed to a diagram of graded vector spaces;
	nothing like this 
	can be hoped to hold in the original diagram of Lie groups.
	}
Since the $G_\a$ are Lie groups,
the indexing partial order is discrete and has minimal elements,
which are now of the form $\H_{S^1}$
or $\H_{G_\a}$ for $G_\a$ simple.
Now an induction is possible. 
For the base case, $Q \H_{G_\a}$ is one-dimensional 
and $s_\a$ is uniquely determined.
For the induction step,
because we have included all quotient groups in the diagram,
each $\H_{G_\a}$ 
is $\smash{\lim_{\b < \a}\H_{G_\b}}$.
As the $s_\b$ have been chosen to make the
squares (\ref{eq:sectionsquare}) commute,
the limit $\smash{s_{\leq \a} \ceq \lim_{\b < \a} s_\b}$ makes sense
and we may take $s_\a$ 
to be the composition 
\[\textstyle\smash{
	Q \H_{G_\a} 
		\isoto 
	\lim_{\b < \a} Q\H_{G_\b} 
		\xmono{s_{\leq \a}} 
	\lim_{\b < \a} \H_{G_\b}
		\isoto 
	\H_{G_\a}.
}\]
This constructs $s_\a$ for all $\a$;
now we may take $\defm s = \lim s_\a$.

It is now clear that (\ref{eq:Cechcolimspace})
is the colimit of the maps (\ref{eq:CartanKapovitch}),
so 
%
%
as colimit is an exact functor, 
we may commute the colimit in (\ref{eq:colimCech}) 
with cohomology
to arrive at an identification of $\Cech(G/K \longinc {}_K G/K)$
with the cohomology of 
the model (\ref{eq:Cechcolimspace})
as claimed.
\epf

From the existence of these pure Sullivan algebras 
we can with little effort extract 
generalizations of results known if $G$ is a Lie group.
The common thread in the proofs 
is that the assumption from the Lie case that the space of generators
is finite-dimensional is actually irrelevant.
We say a \CDGA is \defd{formal} 
if it can be connected through a zig-zag 
of \CDGA quasi-isomorphisms to its own cohomology,
viewed as a \CDGA with differential zero.
A space is formal if its algebra $\APL(X)$
of polynomial differential forms 
is.

\bprop[Cf.~{\cite[pp.~83, 152]{GHVIII}%
\cite[Thm.~12.6.2, p.~211]{onishchik}%
\cite[Thm.~7.4.7,8]{carlsonmonograph}}]\label{thm:Cechformal}
	Let $(G,K)$ be a pair of compact, connected Hausdorff groups.
	The model $(\Cech_K \ox \Cech G,d)$ is formal
	if and only if the ideal of $\Cech_K$ generated
	by the image of 
	$\rho^*\:\wH^{\geq 1}_G \lt \wH^{\geq 1}_K$
	is also generated by a regular sequence 
	contained in this image.
	(For any finitely-generated pure Sullivan algebra
	these conditions are equivalent to the equality 
	$\dim \wh P = \dim P - \dim Q$.)
\eprop

\bprop[Cf.~{\cite[Thms.~A, 3.4]{carlsonfok2018}}]\label{cptisotf}
If a pair $(G,K)$ of compact, connected Hausdorff groups
is isotropy-formal for \v{C}ech cohomology
with rational coefficients, 
then the models of $G/K$ and ${}_K G/K$ considered 
above are formal,
and the cohomology of the latter is isomorphic to
\[
\xt{\wH^*_G}{\wH^*_K\,}{\,\wH^*_K} \,\ox\, \im (\mspace{-1mu}\wH^*_K G_K \lt \Cech G)
\]

\vspace{-1em}

\nd as an $(\Cech_K \ox \Cech_K)$-algebra.
\eprop

\brmk
If $G$ is a Lie group, these propositions mean 
$G/K$ and ${}_K G/K$
are formal in the sense of rational homotopy theory,
but we do not recover this statement in general
because \v{C}ech and singular cohomology $\H X = \H\big(\APL(X))$ 
will differ.
\ermk

\bprop[Cf.~{\cite[p.~218]{cartan1950transgression}\cite[Thm.~2.15.V, p.~73]{GHVIII}\cite[\SS8.4]{onishchik}}]
	Let $(G,K)$ be a pair of compact, connected Hausdorff groups
	such that the model of $G/K$ considered above is formal 
	and let $\widecheck Q < \Cech_G$ be a graded vector subspace 
	sent bijectively by $\rho^*\: \Cech_G \lt \Cech_K$
	to the space spanned by a regular sequence generating 
	$\smash{(\rho^ * \wH^{\geq 1}_G) \lhd \wH_K^*}$.
	Suppose there is a graded subspace $\smash{\widehat Q < \wH_G^*}$,
	meeting $\smash{\widecheck Q}$ trivially,
	such that $\Cech_G$ is the symmetric algebra on 
	$\smash{\widecheck Q \+ \wh Q}$
	and $\rho^*\wh Q \leq 
	(\rho^*\wH^{\geq 1}_G)\.(\rho^ * \wH^{\geq 1}_G)$.
	Then $(G,K)$ is isotropy-formal for \v{C}ech cohomology
	with rational coefficients.
\eprop
\bpf
	We have the liberty to choose the section
	$Q\Cech_G \lt \Cech_G$
	to take $\ker \rho^* + \wH^1_G \.\wH^1_G$
	into $\ker \rho^*$ itself.
	By assumption for each $x$ in a homogeneous basis of $\wh Q$
	we can find $a_j, b_j \in \smash{\wH^1_G}$
	with $x' = x - \sum a_j b_j \in \ker \rho^*$.
	Replacing each $x$ with $x'$, 
	we obtain from $\smash{\widehat Q}$ a different set
	${\widehat Q'}$ such that $\rho^* \wh Q' = 0$
	but $\widecheck Q + \widehat Q'$ still irredundantly generates 
	$\wH^*_K$.
	The suspension maps 
	$\Cech_{G_\a} \epi \wH^{\geq 1}_{G_\a}/ 
	\wH^{\geq 1}_{G_\a}\. 
	\wH^{\geq 1}_{G_\a}
	 \simto P \Cech G_\a \inc \Cech G_\a$
	colimit to a map $\smash{\s\: \Cech_G \lt \Cech G}$
	taking $\smash{\widecheck Q + \wh Q}$ bijectively 
	onto a space of exterior generators $P < \Cech G$.
	If we write $\Cech G = \ext P$
	and $\s \widecheck Q = \cP$ 
	and $\s \wh Q' = \wP$
	then the model $(\Cech_K \ox \Cech_K \ox \Cech G,\td)$
	factors as $\vphantom{X^{X^{X^a}}}(\Cech_K \ox \Cech_K \ox \ext \cP,\td) \ox (\ewP,0)$
	and the resulting map 
	\[
		(\Cech_K \ox \Cech_K \ox \Cech G,\td) \lt
		\Big(\!\mn \quotientbig{\Cech_K \ox \Cech_K\,}{\,(\td \cP)} \ox \ewP, 0\Big)
	\]
	is a quasi-isomorphism.
	Likewise, formality of the model of $G/K$
	implies
		$(\Cech_K \ox \Cech G,d) \lt
		\big({\Cech_K}/{(d \cP)} \ox \ewP, 0\big)$
	is a quasi-isomorphism, so
	the cohomology of the restriction map (\ref{eq:Cechcolimspace})
	modeling the fiber inclusion $G/K \longinc {}_K G/K$
	can be identified with the surjection
	$(\Cech_K \ox \Cech_K)/(\td \cP) \ox \ewP 
	\longepi
	{\Cech_K}/{(d \cP)} \ox \ewP$,
	and $(G,K)$ is isotropy-formal
	for \v{C}ech cohomology. 
\epf

These models give us the desired converse of \Cref{thm:colimboring}
if $K$ is a Lie group. In this case $K_\a \iso K$ far enough up
in the partial order. We can loosen this obvious sufficient 
condition a bit by asking only that the images of the differentials
$\td_\a$ stabilize in a suitable sense.

\bprop\label{thm:stable}
	Let $(G,K)$ be a pair of compact, connected Hausdorff groups,
	presented as a projective limit 
	of compact, connected Lie groups $(G_\a,K_\a)$.
	Endow the Cartan and Kapovitch models with differentials
	such that the obvious ring maps are \DGA homomorphisms,
	as in the proof of \Cref{thm:dumbmodel},
	and suppose there is some index $\w$ such that 
	for all $\a \geq \w$
	the ideal
	$(\td_\a P\H G_\a)$ of $\smash{\H_{K_\a} \ox \H_{K_\a}}$
	is generated by the image of $\td_\w P\H G_\w$
	under
	$(B\phi'_{\a \to \w})^* \ox (B\phi'_{\a \to \w})^*$.
	Then $(G,K)$ is isotropy-formal
	if and only if 
	$(G_\w,K_\w)$ is.
\eprop
\bpf
	This heavy-handed hypothesis ensures that for all $\a \geq \w$
	the primitive elements $\defm \wP^\perp_{\a}$ of $\H G_\a$
	not in the image of $\H G_\w \lt \H G_\a$
	lie in the Samelson subspace 
	for both the Kapovitch and the Cartan algebras,
	so that we can coherently
	tensor-factor 
	the exterior algebra on
	$\defm{\wP^\perp_{\a}} \ceq P\H G_\a / P\H G_\w$,
	equipped with trivial differential,
	out of these models.
	Moreover the induced differentials of the nontrivial
	factors $(\H_{K_\a})^{\ox 2} \ox \H G_\w$
	are determined by the compositions
		\[
		P\H G_\w \lt (\H_{G_\w})^{\ox 2} \longmono
		(\H_{K_\w})^{\ox 2} \longmono 
		(\H_{K_\a})^{\ox 2} \longmono 
		(\H_{K_\b})^{\ox 2},
		\]
		for $\b \geq \a \geq \w$,
	where the last two maps represent each target
	ring as a free module over the source
	so the Kapovitch model 
	$\big((\H_{K_\a})^{\ox 2} \ox \H G_\a,\td_\a\big)$ 
	of ${}_{K_\a} \mn G_\a/{K_\a}$ factors as
	\[
		\xt{\big((\H_{K_\w})^{\ox 2}\big)}{\big((\H_{K_\a})^{\ox 2},0\big)}
			{	\big((\H_{K_\w})^{\ox 2} \ox \H G_\w,
			\td_\w\big)}
		\ox
		(\ext \wP_\a^\perp,0)
	\]
	for $\a \geq \w$ 
	and likewise the Cartan model $(\H_{K_\a} \ox \H G_\a,d_\a)$
	of $G_\a/K_\a$
	factors as
	\[
	\xt{(\H_{K_\w},0)}{(\H_{K_\a},0)}
	{(\H_{K_\w} \ox \H G_\w,d_\w)}
	\ox
	(\ext \wP_\a^\perp,0).
	\]
	In the colimit this describes
	a decomposition of the model of $G/K \longinc {}_K G/K$
	inducing the map
	\quation{\label{lastmap}
		\xt{(\H_{K_\w})^{\ox 2}}{(\Cech_K)^{\ox 2}}
		{\H\big((\H_{K_\w})^{\ox 2} \ox \H G_\w,d_\w\big)}
		\ox
		\ext\wP^\perp
	\lt
		\xt{\H_{K_\w}}{\Cech_K}
		{\H(\H_{K_\w} \ox \H G_\w,d_\w)}
		\ox
		\ext\wP^{\perp}.
	}
	in cohomology.

	If the map
	$\H\big((\H_{K_\w})^{\ox 2} \ox \H G_\w,d_\w\big)
	\lt
	\H(\H_{K_\w} \ox \H G_\w,d_\w)$ 
	of \Cref{thm:dumbmodel}
	arising from $G_\w/K_\w \longinc {}_{K_\w}\mn G_\w/K_\w$
	is surjective, clearly (\ref{lastmap}) is too.
	On the other hand, as 
	$\smash{\Cech_K}$ is a free module $A \ox \H_{K_\w}$ over $\H_{K_\w}$,
	reduction modulo the augmentation ideal of $A$
	makes $\H_{K_\w}$ a $\Cech_K$-module,
	and because $\smash{\Cech_K}$ acts on the ring on the right-hand side
	of (\ref{lastmap}),
	so also does $\smash{(\Cech_K)^{\ox 2}}$ by the reduction
	$\Cech_K \ox \Cech_K \longepi \Q \ox \Cech_K$.
	This action makes (\ref{lastmap}) 
	a map of modules over
	$(\Cech_K)^{\ox 2} \ox \ext\wP^\perp$.
	If this map is surjective,
	then the map obtained by applying 
	$\smash{\big((\H_{K_\w})^{\ox 2} \ox \Q\big) \ox_{(\Cech_K)^{\ox 2}\, \ox\, \ext \wP^\perp} -}$
	is also surjective; but this is just the cohomology of
	the model of $G_\w/K_\w \longinc {}_{K_\w}\mn G_\w/K_\w$
	from \Cref{thm:dumbmodel}.
	\epf

If the $K_\a$ themselves stabilize to $K$, 
then the ideals $(\td_\a P\H G_\a)$ in \Cref{thm:stable} 
must stabilize as well simply since $\HK \ox \HK$ is Noetherian.

\bcor\label{thm:cptoverlie}
Let $(G,K)$ be a pair of compact, connected Hausdorff groups.
If $K$ is a Lie group, 
then $(G,K)$ is \isotf if and only if 
$G$ admits some Lie quotient $\phi\: G \longepi \ol G$ 
such that $K \cap \ker \phi = 1$ and 
$(\ol G, \phi K)$ is \isotf. 
\ecor

This gives us back a version of our circle result.

\bcor\label{thm:cptoverlie}
Let $G$ be a compact, connected Hausdorff group
and $S$ a circle subgroup.
Then $(G,S)$ 
is \isotf if and only if 
$S$ is not contained in the commutator subgroup of $G$
or otherwise there is some Lie quotient $\ol G$ of $G$ 
in which the image of $S$ is a reflected circle, as described in \Cref{algor}.
\ecor

{\footnotesize\bibliography{bibshort} }

\newcommand{\etalchar}[1]{$^{#1}$}
\begin{thebibliography}{BBF{\etalchar{+}}60}

\bibitem[AB84]{AB1984}
Michael~F. Atiyah and Raoul Bott.
\newblock {The moment map and equivariant cohomology}.
\newblock {\em Topology}, 23(1):1--28, 1984.
\newblock \href {https://doi.org/10.1016/0040-9383(84)90021-1}
  {\path{doi:10.1016/0040-9383(84)90021-1}}.

\bibitem[Ada69]{adamsLiebook}
J.~Frank Adams.
\newblock {\em {Lectures on {Lie} groups}}.
\newblock Univ. Chicago Press, 1969.

\bibitem[And62]{andre1962tohoku}
Michel Andr{\'e}.
\newblock {Cohomologie des alg{\`e}bres diff{\'e}rentielles o{\`u} op{\`e}re
  une alg{\`e}bre de {Lie}}.
\newblock {\em T{\^o}hoku Math. J. (2)}, 14(3):263--311, 1962.
\newblock URL:
  \url{http://jstage.jst.go.jp/article/tmj1949/14/3/14_3_263/_article/-char/ja/}.

\bibitem[BBF{\etalchar{+}}60]{borel1960seminar}
Armand Borel, Glen Bredon, Edwin~E. Floyd, Deane Montgomery, and Richard
  Palais.
\newblock {\em {Seminar on transformation groups}}.
\newblock Number~46 in {Ann. of Math. Stud.} Princeton Univ. Press, 1960.
\newblock URL:
  \url{http://indiana.edu/~jfdavis/seminar/Borel,Seminar_on_Transformation_Groups.pdf}.

\bibitem[Bor53]{borelthesis}
Armand Borel.
\newblock {Sur la cohomologie des espaces fibr{\'e}s principaux et des espaces
  homog{\`e}nes de groupes de {L}ie compacts}.
\newblock {\em Ann. of Math. (2)}, 57(1):115--207, 1953.
\newblock URL:
  \url{http://web.math.rochester.edu/people/faculty/doug/otherpapers/Borel-Sur.pdf},
  \href {https://doi.org/10.2307/1969728} {\path{doi:10.2307/1969728}}.

\bibitem[Bor98]{borel1leraybio}
Armand Borel.
\newblock {Jean {Leray} and Algebraic Topology}.
\newblock In Armand Borel, editor, {\em {Jean Leray, selected papers: Oeuvres
  scientifiques}}, volume~1, pages 1--21. Springer and Soc. Math. France, 1998.
\newblock URL:
  \url{http://springer.com/cda/content/document/cda_downloaddocument/9783540609490-c1.pdf}.

\bibitem[Bou68]{bourbakiLie4}
Nicolas Bourbaki.
\newblock {\em {Groupes et alg{\`e}bres de {Lie}, chap. {IV} {\`a} {VI}}}.
\newblock Number~9 in {\'E}l{\'e}ments de math{\'e}matique. Hermann, 1968.

\bibitem[Bou82]{bourbakiLie7}
Nicolas Bourbaki.
\newblock {\em Lie groups and {Lie} algebras: Chapters 7--9}.
\newblock Number~9 in Elements of Mathematics ({\'E}l{\'e}ments de
  math{\'e}matique). Springer (Hermann), 2005 (1975, 1982).
\newblock Transl. by Andrew Pressley.

\bibitem[Bri98]{brion1998eqcohom}
Michel Brion.
\newblock {Equivariant cohomology and equivariant intersection theory}.
\newblock In {\em {Representation theories and algebraic geometry}}, pages
  1--37. Springer, 1998.
\newblock URL:
  \url{http://link.springer.com/chapter/10.1007/978-94-015-9131-7_1}, \href
  {http://arxiv.org/abs/9802063} {\path{arXiv:9802063}}.

\bibitem[Bry17]{MO:Bryant}
Robert Bryant.
\newblock Simple identity on {Lie} algebras in a note of {Koszul}.
\newblock MathOverflow answer, 2017.
\newblock \url{http://mathoverflow.net/q/272042}.

\bibitem[BtD85]{brockertomdieck}
Theodor Br{\"o}cker and Tammo tom Dieck.
\newblock {\em {Representations of compact {Lie} groups}}, volume~98 of {\em
  {Grad. Texts in Math.}}
\newblock Springer, 1985.

\bibitem[BV82]{BV1982}
Nicole Berline and Mich{\`e}le Vergne.
\newblock {Classes caract{\'e}ristiques {\'e}quivariantes. {Formule} de
  localisation en cohomologie {\'e}quivariante}.
\newblock {\em C. R. Acad. Sci. Paris}, 295(2):539--541, Nov 1982.
\newblock \url{http://gallica.bnf.fr/ark:/12148/bpt6k62356694/f77}.

\bibitem[Car51]{cartan1950transgression}
Henri Cartan.
\newblock {La transgression dans un groupe de {Lie} et dans un espace fibr{\'e}
  principal}.
\newblock In {\em {Colloque de topologie (espace fibr{\'e}s), {Bruxelles}
  1950}}, pages 57--71, Li{\`e}ge/Paris, 1951. Centre belge de recherches
  math{\'e}matiques, Georges Thone/Masson et companie.
\newblock \url{http://eudml.org/doc/112227}.

\bibitem[Car72]{carter1972}
Roger~W. Carter.
\newblock {Conjugacy classes in the {Weyl} group}.
\newblock {\em Compos. Math.}, 25(1):1--59, 1972.
\newblock \url{http://numdam.org/item?id=CM_1972__25_1_1_0}.

\bibitem[Car85]{carterbook}
Roger~W. Carter.
\newblock {\em {Finite groups of {Lie} type: Conjugacy classes and complex
  characters}}.
\newblock Wiley, 1985.

\bibitem[Car15]{carlsonmonograph}
Jeffrey~D. Carlson.
\newblock {On the Equivariant Cohomology of Homogeneous Spaces}.
\newblock Manuscript monograph, 2015.
\newblock \url{https://jdkcarlson.github.io/homog_book.pdf}.

\bibitem[CF18]{carlsonfok2018}
Jeffrey~D. Carlson and Chi-Kwong Fok.
\newblock Equivariant formality of isotropy actions.
\newblock {\em J. Lond. Math. Soc.}, Mar 2018.
\newblock \href {http://arxiv.org/abs/1511.06228} {\path{arXiv:1511.06228}},
  \href {https://doi.org/10.1112/jlms.12116} {\path{doi:10.1112/jlms.12116}}.

\bibitem[Che50]{chevalley1950betti}
Claude Chevalley.
\newblock The {Betti} numbers of the exceptional simple {Lie} groups.
\newblock In {\em Proceedings of the International Congress of Mathematicians,
  Cambridge, Mass.}, volume~2, pages 21--24, 1950.

\bibitem[DW98]{dwyer1998Lie}
William~G. Dwyer and Clarence~W. Wilkerson, Jr.
\newblock {The elementary geometric structure of compact {Lie} groups}.
\newblock {\em Bull. London Math. Soc.}, 30(4):337--364, 1998.
\newblock URL:
  \url{http://hopf.math.purdue.edu/Dwyer-Wilkerson/lie/liegroups.pdf}.

\bibitem[DW01]{dwyer2001center}
William~G. Dwyer and Clarence~W. Wilkerson, Jr.
\newblock {Centers and {Coxeter} elements}.
\newblock {\em Contemp. Math.}, 271:53--76, 2001.
\newblock \url{http://www3.nd.edu/~wgd/Dvi/Coxeter.And.Center.pdf}.

\bibitem[FH91]{fultonharris}
William Fulton and Joe Harris.
\newblock {\em {Representation theory}}, volume 129 of {\em {Grad. Texts in
  Math.}}
\newblock Springer, New York, 1991.

\bibitem[FHT01]{FHT}
Yves F{\'e}lix, Steve Halperin, and Jean-Claude Thomas.
\newblock {\em {Rational homotopy theory}}, volume 205 of {\em {Grad. Texts in
  Math.}}
\newblock Springer, 2001.

\bibitem[FOT08]{FOT}
Yves F{\'e}lix, John Oprea, and Daniel Tanr{\'e}.
\newblock {\em {Algebraic models in geometry}}, volume~17 of {\em {Oxford Grad.
  Texts Math.}}
\newblock Oxford Univ. Press, Oxford, 2008.
\newblock \url{www.maths.ed.ac.uk/~v1ranick/papers/tanre.pdf}.

\bibitem[GGK02]{GGK}
Viktor~L. Ginzburg, Victor Guillemin, and Yael Karshon.
\newblock {\em {Moment maps, cobordisms, and {Hamiltonian} group actions}},
  volume~98 of {\em {Math. Surveys Monogr.}}
\newblock Amer. Math. Soc., Providence, RI, 2002.
\newblock URL:
  \url{http://utm.utoronto.ca/~karshony/HUJI/monograph/index-pdf.html}.

\bibitem[GHV76]{GHVIII}
Werner~H. Greub, Stephen Halperin, and Ray Vanstone.
\newblock {\em {Connections, curvature, and cohomology, vol. {III}: Cohomology
  of principal bundles and homogeneous spaces}}.
\newblock Academic Press, 1976.

\bibitem[GKM98]{GKM1998}
Mark Goresky, Robert Kottwitz, and Robert MacPherson.
\newblock {Equivariant cohomology, {Koszul} duality, and the localization
  theorem}.
\newblock {\em Invent. Math.}, 131(1):25--83, 1998.
\newblock \url{http://math.ias.edu/~goresky/pdf/equivariant.jour.pdf},
  \href {https://doi.org/10.1007/s002220050197}
  {\path{doi:10.1007/s002220050197}}.

\bibitem[GN16]{goertschesnoshari2016}
Oliver Goertsches and Sam~Haghshenas Noshari.
\newblock Equivariant formality of isotropy actions on homogeneous spaces
  defined by lie group automorphisms.
\newblock {\em J. Pure Appl. Algebra}, 220(5):2017--2028, 2016.
\newblock \href {http://arxiv.org/abs/1405.2655} {\path{arXiv:1405.2655}},
  \href {https://doi.org/10.1016/j.jpaa.2015.10.013}
  {\path{doi:10.1016/j.jpaa.2015.10.013}}.

\bibitem[Goe12]{goertsches2012isotropy}
Oliver Goertsches.
\newblock {The equivariant cohomology of isotropy actions on symmetric spaces}.
\newblock {\em Doc. Math.}, 17:79--94, 2012.
\newblock \url{emis.ams.org/journals/DMJDMV/vol-17/03.pdf}, \href
  {http://arxiv.org/abs/1009.4079} {\path{arXiv:1009.4079}}.

\bibitem[Hat02]{hatcher}
Allen Hatcher.
\newblock {\em {Algebraic topology}}.
\newblock Cambridge Univ. Press, 2002.
\newblock \url{http://math.cornell.edu/~hatcher/AT/ATpage.html}.

\bibitem[HM06]{hofmannmorriscompact}
Karl~H. Hofmann and Sidney~A. Morris.
\newblock {\em {The structure of compact groups: {A} primer for students --- a
  handbook for the expert}}, volume~25 of {\em {De {Gruyter} Studies in
  Mathematics}}.
\newblock Walter de Gruyter, 2nd revised and augmented edition edition, 2006.

\bibitem[HM07]{hofmannmorrispro}
Karl~H. Hofmann and Sidney~A. Morris.
\newblock {\em {The {Lie} theory of connected pro-{Lie} groups: {A} structure
  theory for pro-{Lie} algebras, pro-{Lie} groups, and connected locally
  compact groups}}, volume~2 of {\em {EMS Tracts Math.}}
\newblock Eur. Math. Soc., 2007.

\bibitem[Hop41]{hopf1941hopf}
Heinz Hopf.
\newblock {{\"U}ber eie {Topologie} der {Gruppen-Mannigfaltigkeiten} und ihre
  {Verallgemeinerungen}}.
\newblock {\em Ann. of Math. (2)}, 42(1):22--52, Jan 1941.
\newblock \href {https://doi.org/10.2307/1968985} {\path{doi:10.2307/1968985}}.

\bibitem[Hsi75]{hsiang}
Wu-Yi Hsiang.
\newblock {\em {Cohomology theory of topological transformation groups}}.
\newblock Springer, 1975.

\bibitem[JK95]{jeffreykirwan1995}
Lisa~C. Jeffrey and Frances~C. Kirwan.
\newblock {Localization for nonabelian group actions}.
\newblock {\em Topology}, 34(2):291--327, 1995.
\newblock \href {http://arxiv.org/abs/alg-geom/9307001}
  {\path{arXiv:alg-geom/9307001}}, \href
  {https://doi.org/10.1016/0040-9383(94)00028-J}
  {\path{doi:10.1016/0040-9383(94)00028-J}}.

\bibitem[Kan01]{kane}
Richard Kane.
\newblock {\em {Reflection groups and invariant theory}}, volume~5 of {\em {C.
  M. S. Books in Mathematics}}.
\newblock Springer, 2001.

\bibitem[Kap04]{kapovitch2002biquotients}
Vitali Kapovitch.
\newblock {A note on rational homotopy of biquotients}.
\newblock 2004.
\newblock \url{http://math.toronto.edu/vtk/biquotient.pdf}.

\bibitem[KN69]{kobayashinomizuII}
Shoshichi Kobayashi and Katsumi Nomizu.
\newblock {\em {Foundations of differential geometry}}, volume~2.
\newblock Interscience Publishers, New York, 1969.

\bibitem[Kos47]{koszul1947homologie}
Jean-Louis Koszul.
\newblock {Sur l'homologie des espaces homog{\`e}nes}.
\newblock {\em C. R. Acad. Sci. Paris}, 225:477--479, Sep 1947.
\newblock \url{http://gallica.bnf.fr/ark:/12148/bpt6k3177x/f477}.

\bibitem[Ler46]{lerayCR1946b}
Jean Leray.
\newblock {Structure de l'anneau d'homologie d'une repr{\'e}sentation}.
\newblock {\em C. R. Acad. Sci. Paris}, 222:1419--1422, Jul 1946.
\newblock \url{http://gallica.bnf.fr/ark:/12148/bpt6k31740/f1419}.

\bibitem[Mil56]{milnor1956universalII}
John~W. Milnor.
\newblock Construction of universal bundles, {II}.
\newblock {\em Ann. of Math. (2)}, 63(3):430--436, May 1956.
\newblock URL:
  \url{http://math.mit.edu/~hrm/18.906/milnor-construction-universal-ii.pdf},
  \href {https://doi.org/10.2307/1970012} {\path{doi:10.2307/1970012}}.

\bibitem[Mil00]{miller2000leray}
Haynes Miller.
\newblock {Leray in {Oflag} {XVIIA}: the origins of sheaf theory, sheaf
  cohomology, and spectral sequences}.
\newblock {\em Gazette des math{\'e}maticiens}, (84 suppl.):17--34, Apr 2000.
\newblock \url{http://www-math.mit.edu/~hrm/papers/ss.pdf}.

\bibitem[Mos53]{mostert1953local}
Paul~S. Mostert.
\newblock Local cross sections in locally compact groups.
\newblock {\em Proc. Amer. Math. Soc.}, 4(4):645--649, 1953.
\newblock \url{http://www.jstor.org/stable/2032540}.

\bibitem[Oni94]{onishchik}
Arkadi~L. Onishchik.
\newblock {\em {Topology of transitive transformation groups}}.
\newblock Johann Ambrosius Barth, 1994.

\bibitem[Sam41]{samelson1941samelson}
Hans Samelson.
\newblock {Beitr{\"a}ge zur {Topologie} der {Gruppen-Mannigfaltigkeiten}}.
\newblock {\em Ann. of Math. (2)}, 42(1):1091--1137, Jan 1941.
\newblock \url{http://jstor.org/stable/1970463}.

\bibitem[Sam49]{samelson3formCR}
Hans Samelson.
\newblock Sur les sous-groupes de dimension 3 des groupes de {Lie} compacts.
\newblock {\em C. R. Acad. Sci. Paris}, 228:630--631, Jan 1949.
\newblock \url{http://gallica.bnf.fr/ark:/12148/bpt6k31801/f630}.

\bibitem[Shi96]{shiga1996equivariant}
Hiroo Shiga.
\newblock {Equivariant de {Rham} cohomology of homogeneous spaces}.
\newblock {\em J. Pure Appl. Algebra}, 106(2):173--183, 1996.
\newblock \href {https://doi.org/10.1016/0022-4049(95)00018-6}
  {\path{doi:10.1016/0022-4049(95)00018-6}}.

\bibitem[Smi67]{smith1967emss}
Larry Smith.
\newblock Homological algebra and the {Eilenberg}--{Moore} spectral sequence.
\newblock {\em Trans. Amer. Math. Soc.}, 129:58--93, 1967.
\newblock \href {https://doi.org/10.2307/1994364} {\path{doi:10.2307/1994364}}.

\bibitem[Spa94]{spanier}
Edwin~H. Spanier.
\newblock {\em Algebraic topology}.
\newblock Springer, 1994.

\bibitem[ST95]{shigatakahashi1995}
Hiroo Shiga and Hideo Takahashi.
\newblock {Remarks on equivariant cohomology of homogeneous spaces}.
\newblock Technical report~17, Tech. Univ. Nagaoka, May 1995.
\newblock \url{dl.ndl.go.jp/info:ndljp/pid/8760355}.

\bibitem[Tay15]{MO:E6}
Jay Taylor.
\newblock {Is this characterization of $(-1)$-eigenspaces of the {Weyl} group
  of {$E_6$} known?}
\newblock MathOverflow answer, 2015.
\newblock \url{http://mathoverflow.net/q/203386}.

\end{thebibliography}

\nd\footnotesize{%
\url{jeffrey.carlson@tufts.edu}
}
\end{document}